\documentclass[a4paper,12pt,leqno]{article}

\usepackage{hyperref}
\usepackage[T1]{fontenc}
\usepackage[english]{babel}
\usepackage{amsmath}
\input xy \xyoption{all}

\newlength\regolo

\newcounter{cnt}[section]
\newcounter{sub}[section]
\newcounter{itm}

\renewcommand\thecnt{\arabic{section}.\arabic{cnt}}
\renewcommand\thesub{\arabic{section}.\arabic{sub}}

\makeatletter

\renewcommand\l@section{\@dottedtocline{1}{0em}{1.5em}}

\newcommand\sezione[2][]{\settowidth\regolo{#1}\ifdim\regolo>0pt\@startsection{section}{1}{0pt}{-\baselineskip}{\baselineskip}{\bfseries\Large}[#1]{#2}\else\@startsection{section}{1}{0pt}{-\baselineskip}{\baselineskip}{\bfseries\Large}[#2]{#2}\fi}

\makeatother

\newenvironment{elenco}{\begin{list}{\roman{itm})}{\setlength\itemindent{0pt}\setlength\labelsep{0.5em}\setlength\labelwidth\leftmargin\addtolength\labelwidth{-\labelsep}\setlength\listparindent{.5\parindent}\setlength\parsep\parskip\setlength\itemsep\medskipamount\setlength\partopsep{0pt}\usecounter{itm}}}{\end{list}}

\newenvironment{equazione}{\addtocounter{sub}{+1}\addtocounter{equation}{+1}\equation\tag{\thesub}}{\endequation}

\newenvironment{equazione+}{\addtocounter{sub}{+1}\addtocounter{equation}{+1}\flalign\tag{\thesub}&\phantom!&&}{\endflalign}

\newenvironment{multiriga}{\addtocounter{sub}{+1}\addtocounter{equation}{+1}\multline\tag{\thesub}}{\endmultline}

\providecommand\qedsymbol{\textit{q.e.d.}}
\newcommand\mathqed{\quad\hbox{\qedsymbol}}
\DeclareRobustCommand\qed{\ifmmode\mathqed\else\leavevmode\unskip\penalty9999\hbox{}\nobreak\hfill\quad\hbox{\qedsymbol}\fi}
\newenvironment{proof}[1][]{\begin{list}{}{\setlength\itemindent\parindent\setlength\labelsep{0pt}\setlength\labelwidth{0pt}\setlength\leftmargin{0pt}\setlength\listparindent\parindent\setlength\parsep\parskip\setlength\partopsep{0pt}}\item\textit{Proof}\settowidth\regolo{#1}\ifdim\regolo>0pt\textit{~~(#1)}\else\fi~~}{\qed\end{list}}

\newenvironment{paragrafo}[1][]{\refstepcounter{cnt}\begin{list}{}{\setlength\itemindent{0pt}\setlength\labelsep\parindent\setlength\labelwidth{0pt}\setlength\leftmargin{0pt}\setlength\listparindent\parindent\setlength\parsep\parskip\setlength\partopsep{0pt}}\item\thecnt.\settowidth\regolo{#1}\ifdim\regolo>0pt\textit{~~#1}\else\fi\ \ }{\end{list}}

\newenvironment{definizione}[1][]{\refstepcounter{cnt}\begin{list}{}{\setlength\itemindent{0pt}\setlength\labelsep{0pt}\setlength\labelwidth{0pt}\setlength\leftmargin{0pt}\setlength\listparindent\parindent\setlength\parsep\parskip\setlength\partopsep{0pt}}\item\textbf{\thecnt~~Definition}~~}{\end{list}}

\newenvironment{lemma}[1][]{\refstepcounter{cnt}\begin{list}{}{\setlength\itemindent{0pt}\setlength\labelsep{0pt}\setlength\labelwidth{0pt}\setlength\leftmargin\parindent\setlength\listparindent\parindent\setlength\parsep\parskip\setlength\partopsep{0pt}}\item\textbf{\thecnt~~Lemma}\settowidth\regolo{#1}\ifdim\regolo>0pt~~(#1)\slshape~~\else\slshape~~\fi}{\end{list}}

\newenvironment{proposizione}[1][]{\refstepcounter{cnt}\begin{list}{}{\setlength\itemindent{0pt}\setlength\labelsep{0pt}\setlength\labelwidth{0pt}\setlength\leftmargin\parindent\setlength\listparindent\parindent\setlength\parsep\parskip\setlength\partopsep{0pt}}\item\textbf{\thecnt~~Proposition}\settowidth\regolo{#1}\ifdim\regolo>0pt~~(#1)\slshape~~\else\slshape~~\fi}{\end{list}}

\newenvironment{theorem}[1][]{\refstepcounter{cnt}\begin{list}{}{\setlength\itemindent{0pt}\setlength\labelsep{0pt}\setlength\labelwidth{0pt}\setlength\leftmargin\parindent\setlength\listparindent\parindent\setlength\parsep\parskip\setlength\partopsep{0pt}}\item\textbf{\thecnt~~Theorem}\settowidth\regolo{#1}\ifdim\regolo>0pt~~(#1)\slshape~~\else\slshape~~\fi}{\end{list}}

\newenvironment{theorem*}[1][]{\begin{list}{}{\setlength\itemindent{0pt}\setlength\labelsep{0pt}\setlength\labelwidth{0pt}\setlength\leftmargin\parindent\setlength\listparindent\parindent\setlength\parsep\parskip\setlength\partopsep{0pt}}\item\textbf{Theorem}\settowidth\regolo{#1}\ifdim\regolo>0pt~~(#1)\slshape~\else\slshape~~\fi}{\end{list}}

\usepackage{amssymb}%
\usepackage{bbold}
\usepackage{mathrsfs}%
\DeclareMathAlphabet{\mathpzc}{T1}{pzc}{m}{it}%
\newcommand\bydef{\stackrel{\text{\tiny\textrm{def}}}{=}}%
\newcommand\gm[1]{\guillemotleft{\tiny~}{#1}{\tiny~}\guillemotright}%
\newcommand\inciso[1]{\nobreakdash---\hspace{0pt}#1\nobreakdash---\hspace{0pt}}%
\newcommand\refcnt[2][]{\settowidth\regolo{#1}\ifdim\regolo>0pt\ifmmode\text{\ref{#1}.\ref{#2}}\else{\ref{#1}.\ref{#2}}\fi\else\ifmmode\text{\ref{#2}}\else{\ref{#2}}\fi\fi}%
\newcommand\refequ[2][]{\settowidth\regolo{#1}\ifdim\regolo>0pt\ifmmode\text{(\ref{#1}.\ref{#2})}\else{(\ref{#1}.\ref{#2})}\fi\else\ifmmode\text{(\ref{#2})}\else{(\ref{#2})}\fi\fi}%
\newcommand\refsez[2][]{\settowidth\regolo{#1}\ifdim\regolo>0pt\S\S{\ref{#1}}\nobreakdash--{\ref{#2}}\else\S{\ref{#2}}\fi}%
\newcommand\Aut{\mathrm{Aut}}%
\newcommand\dimens[1]{\ensuremath{\mathrm{dim}\,#1}}%
\newcommand\End{\mathrm{End}}%
\newcommand\GL{\mathit{GL}}%
\newcommand\Hom{\ensuremath{\mathrm{Hom}}}%
\newcommand\Iso{\mathrm{Iso}}%
\newcommand\kernel[1]{\mathrm{Ker}\,#1}%
\newcommand\C{\mathit C}%
\newcommand\Lebesgue[1][1]{\mathit L^{#1}}%
\newcommand\support[2][]{\ensuremath{\mathrm{supp}_{#1}\,#2}}%
\newcommand\T[2][]{\ensuremath{\mathit T_{#1}\mspace{1.5mu}#2}}%
\newcommand\test[1][\infty]{{\mathit C_{\mathit c}}^{#1}}
\newcommand\epito\twoheadrightarrow%
\newcommand\from\leftarrow%
\newcommand\infrom\hookleftarrow%
\newcommand\into\hookrightarrow%
\newcommand\isoto{\stackrel\thicksim\to}%
\newcommand\longto\longrightarrow%
\newcommand\xfrom[1]{\xleftarrow{#1}}%
\newcommand\xto[1]{\xrightarrow{#1}}%

\newcommand\Id{\mathit{Id}}
\newcommand\Ob{\mathrm{Ob}}

\newcommand\VectorSpaces{\ensuremath{\left\{\mathrm{vector\:spaces}\right\}}}

\newcommand\ComplexVectorSpaces{\ensuremath{\left\{\mathrm{complex\:vector\:spaces}\right\}}}

\newcommand\SheavesOfModules[1]{\ensuremath{\left\{\mathrm{sheaves\,of\:{#1}\text-modules}\right\}}}

\newcommand\aeq\Leftrightarrow%
\newcommand\seq\Rightarrow%
\newcommand\pt{\ensuremath{\star}}%
\newcommand\id{\mathit{id}}%
\newcommand\pr{\mathit{pr}}%
\newcommand\nR{\ensuremath{\mathbb{R}}}%
\newcommand\nC{\ensuremath{\mathbb{C}}}%
\newcommand\displaycap[3]{\ensuremath{\overset{#2}{\underset{#1}{\displaystyle\bigcap}}\,#3}}%
\newcommand\displaysum[3]{\ensuremath{\overset{#2}{\underset{#1}{\displaystyle\sum}}\,#3}}%
\newcommand\modulo[1]{\ensuremath{\left|#1\right|}}%
\newcommand\bigmod[1]{\ensuremath{\bigl|#1\bigr|}}%
\newcommand\scalare[2]{\ensuremath{\left\langle#1,#2\right\rangle}}%
\newcommand\bigsca[2]{\ensuremath{\bigl\langle#1,#2\bigr\rangle}}%
\newcommand\parteRe[1]{\ensuremath{\Re\mathfrak e#1}}%
\newcommand\txtcap[3]{\ensuremath{\overset{#2}{\underset{#1}{\textstyle\cap}}\,#3}}%
\newcommand\txtcup[3]{\ensuremath{\overset{#2}{\underset{#1}{\textstyle\cup}}\,#3}}%
\newcommand\txtcoprod[3]{\ensuremath{\overset{#2}{\underset{#1}{\textstyle\coprod}}\,#3}}%
\newcommand\txtsum[3]{\ensuremath{\overset{#2}{\underset{#1}{\textstyle\sum}}\,#3}}%
\newcommand\can\cong%
\newcommand\iso\approx%

\newcommand\sheafconst[2][]{\settowidth{\regolo}{\ensuremath{#1}}\ifdim\regolo>0pt\ensuremath{\underline{#2}_{#1}}\else\ensuremath{\underline{#2}}\fi}

\newcommand\smooth[1][]{\settowidth{\regolo}{\ensuremath{#1}}\ifdim\regolo>0pt\ensuremath{\mathscr C^\infty_{#1}}\else\ensuremath{\mathscr C^\infty}\fi}

\newcommand\G{\ensuremath{\mathcal G}}
\renewcommand\H{\ensuremath{\mathcal H}}
\newcommand\K{\ensuremath{\mathcal K}}

\newcommand\s[1][]{\settowidth{\regolo}{\ensuremath{#1}}\ifdim\regolo>0pt\ensuremath{\mathit s\mspace{.8mu}#1}\else\ensuremath{\mathit s}\fi}

\renewcommand\t[1][]{\settowidth{\regolo}{\ensuremath{#1}}\ifdim\regolo>0pt\ensuremath{\mathit t\mspace{.8mu}#1}\else\ensuremath{\mathit t}\fi}

\renewcommand\c{\ensuremath{\mathit c}}
\newcommand\p{\ensuremath{\mathit p}}
\renewcommand\u{\ensuremath{\mathit u}}
\newcommand\imap{\ensuremath{\mathit i}}

\newcommand\mca[2][]{\settowidth{\regolo}{\ensuremath{#1}}\ifdim\regolo>0pt\ensuremath{#2^{\scriptscriptstyle(#1)}}\else\ensuremath{#2{_{\s}\times_{\t}}#2}\fi}

\newcommand\Kt[1][C]{\ensuremath{\mathcal{#1}}}
\newcommand\TU{\ensuremath{\mathsf1}}

\newcommand\stack[2][]{\settowidth{\regolo}{\ensuremath{#1}}\ifdim\regolo>0pt\ensuremath{\mathfrak{#2}(#1)}\else\ensuremath{\mathfrak{#2}}\fi}

\newcommand\xstack[2][]{\settowidth{\regolo}{\ensuremath{#1}}\ifdim\regolo>0pt\ensuremath{\mathfrak{#2}\mspace{-1.4mu}\left(#1\right)}\else\ensuremath{\mathfrak{#2}}\fi}

\newcommand\sheafhom[4][]{\ensuremath{\mathscr H\mspace{-5.4mu}\mathit{om}^{\stack{#1}}_{#2}(#3,#4)}}

\newcommand\sections[2][]{\settowidth{\regolo}{\ensuremath{#1}}\ifdim\regolo>0pt\ensuremath{\boldsymbol\Gamma_{#1}#2}\else\ensuremath{\boldsymbol\Gamma#2}\fi}

\newcommand\Des[2][]{\settowidth{\regolo}{\ensuremath{#1}}\ifdim\regolo>0pt\ensuremath{\underline{\mathcal D\mspace{-2.6mu}\mathit{es}}^{\mathfrak{#1}}(#2)}\else\ensuremath{\underline{\mathcal D\mspace{-2.6mu}\mathit{es}}(#2)}\fi}

\newcommand\V[2][]{\settowidth{\regolo}{\ensuremath{#1}}\ifdim\regolo>0pt\ensuremath{\mathit V^{\mathfrak{#1}}(#2)}\else\ensuremath{\mathit V(#2)}\fi}

\newcommand\SpV{\ensuremath{\underline{\mathcal V\mspace{-3.0mu}\mathit{ec}}}}

\newcommand\VB[2][\nC]{\ensuremath{{\SpV_{#1}}^\infty(#2)}}

\newcommand\ReV{\ensuremath{\underline{\phantom{\smash[b]{\mathcal R\mspace{-3.0mu}\mathit{ep}}}}\mspace{-30.55mu}{\mathcal R\mspace{-3.0mu}\mathit{ep}}}}

\newcommand\R[2][]{\settowidth{\regolo}{\ensuremath{#1}}\ifdim\regolo>0pt\ensuremath{\mathit R^{\mathfrak{#1}}(#2)}\else\ensuremath{\mathit R(#2)}\fi}

\newcommand\forgetfulfunctor[2][]{\settowidth\regolo{\ensuremath{#1}}\ifdim\regolo>0pt\ensuremath{\mathsf F^{\mathfrak{#1}}_{#2}}\else\ensuremath{\mathsf F_{#2}}\fi}

\newcommand\Av[1]{\mathrm{Av}_{#1}}

\newcommand\envelopemorph[2][]{\settowidth\regolo{\ensuremath{#1}}\ifdim\regolo>0pt\ensuremath{\boldsymbol\pi^{\mathfrak{#1}}_{#2}}\else\ensuremath{\boldsymbol\pi_{#2}}\fi}

\newcommand\f[2][]{\settowidth{\regolo}{\ensuremath{#1}}\ifdim\regolo=0pt\ensuremath{\sections{\mathscr{#2}}}\else\ensuremath{\mathscr{#2}}\fi}

\newcommand\E{\f[*]E}

\newcommand\Euc{\underline{\mathcal E\mspace{-3.0mu}\mathit{uc}}}
\newcommand\Hilb{\underline{\mathcal H\mspace{-2.2mu}\mathit{ilb}}}

\newcommand\EB[2][\infty]{\ensuremath{{\Euc}^{#1}(#2)}}
\newcommand\HB[2][\infty]{\ensuremath{{\Hilb}^{#1}(#2)}}

\newcommand\fifu[1][\omega]{\ensuremath{\boldsymbol{#1}}}

\newcommand\tannakian[2][]{\settowidth{\regolo}{\ensuremath{#1}}\ifdim\regolo=0pt\settowidth{\regolo}{\ensuremath{#2}}\ifdim\regolo>0pt\ensuremath{\mathcal T(#2)}\else\ensuremath{\mathcal T}\fi\else\settowidth{\regolo}{\ensuremath{#2}}\ifdim\regolo>0pt\ensuremath{\mathcal T^{\mathfrak{#1}}(#2)}\else\ensuremath{\mathcal T^{\mathfrak{#1}}}\fi\fi}

\newcommand\Rf{\ensuremath{\mathscr R}}

\begin{document}

\title{A Tannaka Theorem for Proper Lie Groupoids}
\author{Giorgio Trentinaglia\protect\footnote{During the preparation of this paper, the author was partially supported by a grant of the foundation ``Fon\-da\-zio\-ne \mbox{Ing.~Aldo~Gini}''.}}
\date{}
\maketitle

\begin{abstract}
\noindent By replacing the category of smooth vector bundles over a manifold with the category of what we call smooth Euclidean fields, which is a proper enlargement of the former, and by considering smooth actions of Lie groupoids on smooth Euclidean fields, we are able to prove a Tannaka duality theorem for proper Lie groupoids. The notion of smooth Euclidean field we introduce here is the smooth, finite dimensional analogue of the usual notion of continuous Hilbert field.

\end{abstract}

\section*{Introduction}
Classical \mbox{Tannaka}--\mbox{Kre\v{\i}n} duality theory leads to the result that a compact group can be reconstructed from a purely discrete, algebraic object, namely the ring of its continuous finite dimensional representations or, more precisely, the algebra of its representative functions. Compare \cite{BtD85}. The same theory can be efficiently recast in categorical terms. This alternative point of view on Tannaka duality stems from \mbox{Grothendieck's} theory of motives in algebraic geometry \cite{Saa72,DeMi82,De90}. In this approach one starts by considering, for an arbitrary locally compact group $G$, the category formed by the continuous representations of $G$ on finite dimensional vector spaces, endowed with the symmetric monoidal structure arising from the usual tensor product of representations, and then one tries to recover $G$ as the group of all tensor preserving natural endomorphisms of the standard forgetful functor which assigns each $G$-module the underlying vector space. See for instance \cite{JoSt91}. When $G$ is a compact Lie group, in particular, it follows that $G$ can be reconstructed in this way up to isomorphism, as the $\C^\infty$ manifold structure of a Lie group is determined by the underlying topology.

It is natural to ask for a generalization of the aforesaid duality theory to the realm of Lie groupoids, in which proper groupoids are expected to play the same role as compact groups. When trying to extend Tannaka theory from Lie groups to Lie groupoids, however, one is first of all confronted with the problem of choosing a suitable notion of representation for the latter. Now, the notion of smooth or, equivalently, continuous finite dimensional representation has an obvious naive generalization to the Lie groupoid setting: the representations of a Lie groupoid \G\ could be defined to be the Lie groupoid homomorphisms $\G \to \GL(E)$ of \G\ into the linear groupoids associated with smooth vector bundles over its base manifold. Unfortunately, this naive generalization turns out to be inadequate for the purpose of carrying forward Tannaka duality to Lie groupoids, cfr \cite{2008,LuOl01,Henr05}. This state of affairs forces us to introduce a different notion of representation for Lie groupoids. It seems reasonable to require that the new notion should be as close as possible to the naive notion recalled above, and that moreover in the case of groups one should recover the usual notion of smooth representation on a finite dimensional vector space.

In this paper we work out the problems raised in the preceding paragraph. To begin with, we construct, for each smooth manifold $X$, a category whose objects we call smooth Euclidean fields over $X$. Our notion of smooth Euclidean field is the analogue, in the smooth and finite dimensional setting we confine ourselves to, of the notion of continuous Hilbert field introduced by \mbox{Dixmier} and \mbox{Douady} \cite{DixDou63}. The category of smooth Euclidean fields over $X$ is, for every paracompact manifold $X$, a proper enlargement of the category of smooth vector bundles over $X$. One can straightforwardly define a notion of representation of a Lie groupoid on a smooth Euclidean field; such representations form, for each Lie groupoid \G, a symmetric monoidal category which is connected to the category of smooth Euclidean fields over the base manifold of \G\ by a canonical forgetful functor. From this functor we obtain, by generalizing the construction mentioned at the beginning, a groupoid. This ``reconstructed'' groupoid\inciso{the Tannakian groupoid of \G, as we call it}comes equipped with a natural candidate for a differentiable structure on its space of arrows, namely a sheaf of algebras of continuous real valued functions stable under composition with arbitrary smooth functions of several variables. A space endowed with such a structure constitutes what we call a $\C^\infty$-space. There is a canonical homomorphism of \G\ into its Tannakian groupoid, which proves to be also a morphism of $\C^\infty$-spaces. Now, our duality result (Theorem \ref{xx.18}) can be stated as follows:
\begin{theorem*}
For a proper Lie groupoid \G, the canonical homomorphism of \G\ into its Tannakian groupoid is an isomorphism of both groupoids and $\C^\infty$-spaces. It follows that the Tannakian groupoid itself is a Lie groupoid, isomorphic to \G.
\end{theorem*}
Our argument is complementary to the proof of the classical Tannaka duality theorem. Most efforts are directed towards showing how the classical theorem implies the surjectivity of the above-mentioned canonical isomorphism and then towards establishing the claim about the $\C^\infty$-space structure. By contrast, the fact that the canonical homomorphism is injective is a direct consequence of a theorem of \mbox{N.~T.~Zung} \cite{Zu06}; \mbox{Zung's} theorem may in fact be regarded as a \mbox{Peter}--\mbox{Weyl} theorem for proper Lie groupoids. Compare \cite{2008}.

Many of the reasonings leading to our duality theorem, although of course not all of them, also apply to the representations of proper Lie groupoids on vector bundles. Since from the very beginning of our research we were equally interested in studying such representations, we found it convenient to provide a general theoretical framework where the diverse approaches to the representation theory of Lie groupoids could take their appropriate place, so as to state our results in a uniform language. The outcome of this demand was the theory of smooth tensor stacks. Smooth vector bundles and smooth Euclidean fields are two examples of a smooth tensor stack. Each smooth tensor stack gives rise to a corresponding notion of representation for Lie groupoids, and then for each Lie groupoid one obtains, by the same general procedure outlined above, a corresponding Tannakian groupoid. What this groupoid looks like will depend very much, in general, on the initial choice of a smooth tensor stack, as we pointed out already in the course of this introduction.

\bigskip

\noindent\textsc{Acknowledgements.} The problem of proving a Tannaka duality theorem for proper Lie groupoids was originally suggested to the author by \mbox{I.~Moerdijk}, who also made several useful remarks on an earlier draft of this paper. Besides, the author would like to thank \mbox{M.~Crainic} and \mbox{N.~T.~Zung} for helpful discussions.

\tableofcontents

\sezione{Proper Lie Groupoids}\label{i}
The present section is introductory. Its purpose is to recall some background notions and to fix some notation that we will be using throughout the paper. The reader is advised to consult \cite{MoeMrc03,Bre72,Cra03,Zu06} for reference; other sources include \cite{MoeRey91} and \cite{CraMoe01}.

\bigskip

\noindent The term {groupoid} refers to a small category where every arrow is invertible. A {Lie groupoid} can be approximately described as an internal groupoid in the category of smooth manifolds. To construct a Lie groupoid \G\ one has to give a pair of manifolds of class $\C^\infty$ \mca[0]\G\ and \mca[1]\G, respectively called {manifold of objects} and {manifold of arrows,} and a list of smooth maps called {structure maps.} The basic items in this list are the {source} map $\s: \mca[1]\G \to \mca[0]\G$ and the {target} map $\t: \mca[1]\G \to \mca[0]\G$, which have to meet the requirement that the fibred product $\mca[2]\G = \mca{\mca[1]\G}$ exists in the category of $\C^\infty$ manifolds. Then one has to give a {composition} map $\c: \mca[2]\G \to \mca[1]\G$, a {unit} map $\u: \mca[0]\G \to \mca[1]\G$ and an {inverse} map $\imap: \mca[1]\G \to \mca[1]\G$, for which the familiar algebraic laws must be satisfied.

\textit{Terminology and Notation:} The points $x = \s(g)$ and $x' = \t(g)$ are resp.\ called the {source} and the {target} of the {arrow} $g$. We let $\G(x,x')$ denote the set of all the arrows whose source is $x$ and whose target is $x'$, and we use the abbreviation $\G|_x$ for the {isotropy} or {vertex} group $\G(x,x)$. Notationally, we will often identify a point $x \in \mca[0]\G$ and the corresponding unit arrow $\u(x) \in \mca[1]\G$. It is costumary to write ${g'\cdot g}$ or $g'g$ for the composition $\c(g',g)$ and $g^{-1}$ for the inverse $\imap(g)$.

Our description of the notion of Lie groupoid is still incomplete. It turns out that a couple of additional requirements are needed in order to get a reasonable definition.

Recall that a manifold $M$ is said to be {paracompact} if it is Hausdorff and there exists an ascending sequence of open subsets with compact closure $\cdots \subset U_i \subset \overline U_i \subset U_{i+1} \subset \cdots$ such that $M = \txtcup{i=0}\infty{U_i}$. A Hausdorff manifold is paracompact if and only if it possesses a countable basis of open subsets. Any open cover of a paracompact manifold admits a locally finite refinement. Any paracompact manifold admits partitions of unity of class $\C^\infty$ (subordinated to open covers).

In order to make the fibred product \mca{\mca[1]\G} meaningful as a manifold and for other purposes related to our study, we shall include the following additional conditions in the definition of Lie groupoid:
\begin{elenco}
\item[\textsl{1.}]The source map $\s: \mca[1]\G \to \mca[0]\G$ is a submersion with Hausdorff fibres;
\item[\textsl{2.}]The manifold \mca[0]\G\ is paracompact.
\end{elenco}
Note that we do not require the manifold of arrows \mca[1]\G\ to be Hausdorff or paracompact. Actually, this manifold is neither Hausdorff nor second countable in many examples of interest. The first condition implies at once that the domain of the composition map is a submanifold of the Cartesian product ${\mca[1]\G \times \mca[1]\G}$ and that the target map is a submersion with Hausdorff fibres. Thus, the source fibres $\G(x,\text-) = \s^{-1}(x)$ and the target fibres $\G(\text-,x') = \t^{-1}(x')$ are closed Hausdorff submanifolds of \mca[1]\G. A Lie groupoid \G\ is said to be {Hausdorff} if the manifold \mca[1]\G\ is Hausdorff.

\textit{Some more Terminology:} The manifold \mca[0]\G\ is usually called the {base} of the groupoid \G. One also says that \G\ is a groupoid over the manifold \mca[0]\G. We shall often use the notation $\G^x = \G(x,\text-) = \s^{-1}(x)$ for the fibre of the source map over a point $x \in \mca[0]\G$. More generally, we shall write
\begin{equazione}\label{i1}
\G(S,S') = \bigl\{g \in \mca[1]\G: \s(g) \in S\text{~\&~} \t(g) \in S'\bigr\}\text, \quad \G|_S = \G(S,S)
\end{equazione}
and $\G^S = \G(S,\text-) = \G(S,\mca[0]\G) = \s^{-1}(S)$ for all subsets $S, S' \subset \mca[0]\G$.

\textit{Example:} Let $G$ be a Lie group acting smoothly (from the left) on a manifold $M$. Then the {action} (or {translation}) {groupoid ${G\ltimes M}$} is defined to be the Lie groupoid over $M$ whose manifold of arrows is the Cartesian product ${G\times M}$, whose source and target maps are respectively the projection onto the second factor $(g,x) \mapsto x$ and the action $(g,x) \mapsto {g x}$, and whose composition law is the operation
\begin{equazione}\label{i6}
{(g',x')(g,x)} = (g'g,x)\text.
\end{equazione}
There is a similar construction ${M\rtimes G}$ associated with right actions.

A {homomorphism of Lie groupoids} is a smooth functor. More precisely, a homomorphism $\varphi: \G \to \H$ consists of two smooth maps $\mca[0]\varphi: \mca[0]\G \to \mca[0]\H$ and $\mca[1]\varphi: \mca[1]\G \to \mca[1]\H$ compatible with the groupoid structure in the sense that ${\s\circ \mca[1]\varphi} = {\mca[0]\varphi \circ \s}$, ${\t\circ \mca[1]\varphi} = {\mca[0]\varphi \circ \t}$ and $\mca[1]\varphi({g'\cdot g}) = {\mca[1]\varphi(g') \cdot \mca[1]\varphi(g)}$. Lie groupoids and their homomorphisms form a category. A homomorphism $\varphi: \G \to \H$ is said to be a {Morita equivalence} when
\begin{equazione}\label{ME1}
\begin{split}
\xymatrix@C=45pt{\mca[1]\G\ar[d]_{(\s,\t)}\ar[r]^{\mca[1]\varphi} & \mca[1]\H\ar[d]^{(\s,\t)} \\ {\mca[0]\G \times \mca[0]\G}\ar[r]^{\mca[0]\varphi \times \mca[0]\varphi} & {\mca[0]\H \times \mca[0]\H}}
\end{split}
\end{equazione}
is a pullback diagram in the category of $\C^\infty$ manifolds and the map
\begin{equazione}\label{ME2}
{\t_\H \circ \pr_2}: {\mca[0]\G {_{\mca[0]\varphi}\times_{\s_\H}} \mca[1]\H} \to \mca[0]\H
\end{equazione}
is a surjective submersion.

There is also a notion of {topological groupoid.} This is just an internal groupoid in the category of topological spaces and continuous mappings. In the continuous case the definition is much simpler, since one need not worry about the domain of definition of the composition map. With the obvious notion of homomorphism, topological groupoids constitute a category.

Let \G\ be a Lie groupoid and let $x$ be a point of its base manifold \mca[0]\G. The {orbit of \G} (or {\G-orbit}) {through $x$} is the subset
\begin{equazione}\label{i7}
{\G x} \bydef {\G\cdot x} \bydef \t\bigl(\G^x\bigr) = \{x' \in \mca[0]\G| \exists g: x \to x'\}\text.
\end{equazione}
Note that the isotropy group $\G|_x$ acts from the the right on the manifold $\G^x$. This action is clearly free and transitive along the fibres of the restriction of the target map to $\G^x$. The following facts hold, cfr \cite{MoeMrc03} p.~115: \textsl{(a)}~$\G(x,x')$ is a closed submanifold of \mca[1]\G\ \textsl{(b)}~$G_x = \G|_x$ is a Lie group \textsl{(c)}~the \G-orbit through $x$ is an immersed submanifold of \mca[0]{\G} and the target map $\t: \G^x \to {\G x}$ determines a principal $G_x$-bundle over it (the set ${\G x}$ can obviously be identified with the homogeneous space ${\G^x/G_x}$, and it can be proved that there exists a possibly non-Hausdorff manifold structure on this quotient space such that the quotient map turns out to be a principal bundle).

\begin{paragrafo}[$\C^\infty$-Spaces.]\label{C^infty-Sp}
Recall that a functionally structured space is a topological space $X$ endowed with a sheaf $\mathscr F$ of real algebras of continuous real valued functions on $X$ (functional structure). Compare for instance \cite{Bre72}, p.~297. There is an obvious notion of morphism for such spaces.\footnote{Algebraic geometers would say that a morphism of functionally structured spaces is a (continuous) mapping inducing a morphism of (locally) ringed spaces.}

Let $\mathscr F$ be an arbitrary functional structure on a topological space $X$. We shall let $\mathscr F^\infty$ denote the sheaf of continuous real valued functions on $X$ generated by the presheaf
\begin{multiriga}\label{C^infty-Sp1}
U\: \mapsto\: \bigl\{f(a_1|_U,\ldots,a_d|_U): \text{~$f: \nR^d \to \nR$ of class $\C^\infty$,}\\
a_1, \ldots, a_d \in \mathscr F(U)\bigr\}\text.
\end{multiriga}
Here the expression $f(a_1|_U,\ldots,a_d|_U)$ stands of course for the function $u \mapsto f\bigl(a_1(u),\ldots,a_d(u)\bigr)$ on $U$. The pair $(X,\mathscr F^\infty)$ constitutes a functionally structured space to which we shall refer as a \textit{$\C^\infty$ functionally structured space} or, in short, \textit{$\C^\infty$-space.} More correctly, a $\C^\infty$-space is a functionally structured space $(X,\mathscr F)$ such that $\mathscr F = \mathscr F^\infty$.\footnote{A more general notion of ``{$\C^\infty$-ring}'' was introduced by \mbox{Moerdijk} and \mbox{Reyes} in the context of smooth infinitesimal analysis \cite{MoeRey86,MoeRey91}. What we are considering here is a special instance of that notion, namely a $\C^\infty$-ring of continuous functions on a topological space. For simplicity, we choose to work within the subcategory of such $\C^\infty$-rings.} Observe that smooth manifolds can be defined as topological spaces endowed with a $\C^\infty$ functional structure locally isomorphic to that of smooth functions on $\nR^n$. $\C^\infty$-Spaces have, in general, better categorical properties than smooth manifolds. Note that the latter form, within $\C^\infty$-spaces, a full subcategory.
\end{paragrafo}

\begin{paragrafo}[$\C^\infty$-Groupoids.]\label{C^infty-Gpd}
Let us start by observing that if $(X,\mathscr F)$ is a $\C^\infty$-space then so is $(S,\mathscr F|_S)$ for any subspace $S$ of $X$, where $\mathscr F|_S = {{i_S}^*\mathscr F}$ denotes the inverse image of $\mathscr F$ along the inclusion $i_S: S \into X$. [Recall that, for an arbitrary continuous mapping $f: S \to T$ into a functionally structured space $(T,\mathscr T)$, ${f^*\mathscr T}$ denotes the functional sheaf on $S$ formed by the functions which are locally the pullback along $f$ of functions in $\mathscr T$.]

We note next that if $(X,\mathscr F)$ and $(Y,\mathscr G)$ are two functionally structured spaces then so is their Cartesian product endowed with the sheaf ${\mathscr F\otimes\mathscr G}$ locally generated by the functions $({\varphi\otimes\psi})(x,y) = {\varphi(x) \psi(y)}$. It follows that $({\mathscr F^\infty\otimes\mathscr G^\infty})^\infty$ is a $\C^\infty$ functional structure on ${X\times Y}$ turning this into the product of $(X,\mathscr F^\infty)$ and $(Y,\mathscr G^\infty)$ in the category of $\C^\infty$-spaces. The preceding considerations imply that the category of $\C^\infty$-spaces is closed under fibred products (pullbacks). Notice that when $X$ and $Y$ are smooth manifolds and $S \subset X$ is a submanifold one recovers the correct manifold structures, so that all these constructions for $\C^\infty$-spaces agree with the usual ones on manifolds whenever the latter make sense.

We shall use the term \textit{$\C^\infty$-groupoid} to indicate a groupoid whose sets of objects and arrows are each endowed with the structure of a $\C^\infty$-space so that all the maps arising from the groupoid structure (source, target, composition, unit section, inverse) are morphisms of $\C^\infty$-spaces. The base space $X$ is always a smooth manifold in practice, with $\C^\infty$ functional structure given by the sheaf of smooth functions on $X$. Every Lie groupoid is, in particular, an example of a $\C^\infty$-groupoid.
\end{paragrafo}

A Lie (or topological or $\C^\infty$) groupoid \G\ is said to be \textit{proper} if \G\ is Hausdorff and the combined source--target map $(\s,\t): \mca[1]\G \to {\mca[0]\G \times \mca[0]\G}$ is proper (in the familiar sense: the inverse image of a compact subset is compact). When \G\ is a proper Lie groupoid over a manifold $M$, every \G-orbit is in fact a closed submanifold of $M$.

Normalized Haar systems on proper Lie groupoids are the analogue of Haar probability measures on compact Lie groups. We now recall the definition and the construction of Haar systems on proper Lie groupoids. Our exposition is based on \cite{Cra03}. Let \G\ be a Lie groupoid over a manifold $M$.

\begin{definizione}\label{PHaarSyst}
A \textit{positive Haar system} on \G\ is a family of positive Radon measures $\{\mu^x\}$ ($x \in M$), each one supported by the respective source fibre $\G^x = \G(x,\text-) = \s^{-1}(x)$, satisfying the following conditions:
\begin{elenco}
\item${\int \varphi\, \mathit d\mu^x} > 0$ for all nonnegative real $\varphi \in \test(\G^x)$ with $\varphi \neq 0$;
\item for each $\varphi \in \test(\mca[1]\G)$ the function $\Phi$ on $M$ defined by setting
\begin{equazione}\label{PHaarSyst2}
\Phi(x) = {\int_{\G^x} \varphi|_{\G^x}\, \mathnormal d\mu^x}
\end{equazione}
is of class $\C^\infty$;
\item(right invariance) for all $g \in \G(x,y)$ and $\varphi \in \test(\G^x)$ one has
\begin{equazione}\label{PHaarSyst3}
{\int_{\G^y} {\varphi\circ\tau^g}\, \mathit d\mu^y} = {\int_{\G^x} \varphi\, \mathit d\mu^x}
\end{equazione}
where $\tau^g: \G(y,\text-) \to \G(x,\text-)$ denotes right translation $h \mapsto hg$.
\end{elenco}
\end{definizione}

The existence of positive Haar systems on a Lie groupoid \G\ can be established when \G\ is {proper.} One way to do this is the following. One starts by fixing a Riemann metric on the vector bundle $\mathfrak g \to M$, where $\mathfrak g$ is the Lie algebroid of \G\ (cfr \cite{Cra03} or Chapter 6 of \cite{MoeMrc03}; note the use of paracompactness). Right translations determine isomorphisms $\T{\G(x,\text-)} \iso {\t^*\mathfrak g}|_{\G(x,\text-)}$ for all $x \in M$. These isomorphisms can be used to induce, on the source fibres $\G(x,\text-)$, Riemann metrics whose associated volume densities provide the desired system of measures.

\begin{definizione}\label{NHaarSyst}
A \textit{normalized Haar system} on \G\ is a family of positive Radon measures $\{\mu^x\}$ ($x \in M$), each one with support concentrated in the respective source fibre $\G^x$, enjoying the following properties: \textsl{(a)}~All smooth functions on $\G^x$ are integrable with respect to $\mu^x$, that is to say
\begin{equazione}
\C^\infty(\G^x) \subset \Lebesgue(\mu^x)
\end{equazione}
\textsl{(b)}~Condition \textsl{ii)}, respectively \textsl{iii)} of the preceding definition holds for an arbitrary smooth function $\varphi$ on \mca[1]{\G}, respectively $\G^x$ \textsl{(c)}~The following normalization condition is satisfied:
\begin{elenco}
\item[\textsl{i*)}]${\int \mathit d\mu^x} = 1$ for every $x \in M$.
\end{elenco}
\end{definizione}

Every proper Lie groupoid admits normalized Haar systems. One can prove this by using a cut-off function, namely a positive, smooth function $c$ on the base $M$ of the groupoid such that the source map restricts to a proper map on $\support{c\circ\t}$ and ${\int {c\circ\t}\, \mathit d\nu^x} = 1$ for all $x \in M$, where $\{\nu^x\}$ is an arbitrary positive Haar system chosen in advance. The system of positive measures $\mu^x \equiv {({c\circ\t})\, \nu^x}$ will have the desired properties.

\begin{paragrafo}[Zung's theorem.]\label{ZungThm}
Let \G\ be a Lie groupoid and let $M$ be its base manifold. We say that a submanifold $N$ of $M$ is a {slice at a point $z \in N$} if the orbit immersion ${\G z} \into M$ is transversal to $N$ at $z$. A submanifold $S$ of $M$ will be called a {slice} if it is a slice at all of its points. Note that if $N$ is a submanifold of $M$ and $g \in \G^N = \s^{-1}(N)$ then $N$ is a slice at $z = \s(g)$ if and only if the intersection ${\G^N \cap \t^{-1}(z')}, z' = \t(g)$ is transversal at $g$. From this remark it follows that for each submanifold $N$ the subset of all points at which $N$ is a slice forms an open subset of $N$. If a submanifold $S$ of $M$ is a slice then the intersection ${\s^{-1}(S) \cap \t^{-1}(S)}$ is transversal, so that the restriction $\G|_S$ is a Lie groupoid over $S$; moreover, ${\G\cdot S}$ is an invariant open subset of $M$. For the proof of the following result, we refer the reader to \cite{Zu06}.
\begin{theorem*}[N.~T.~Zung]
Let \G\ be a proper Lie groupoid. Let $x$ be a base point which is not moved by the tautological action of $\G$ on its own base. Then there exists a continuous linear representation $G \to \GL(V)$ of the isotropy group $G = \G|_x$ on a finite dimensional vector space $V$ such that for some open neighbourhood $U$ of $x$ one can find an isomorphism of Lie groupoids $\G|_U \iso {G \ltimes V}$ which makes $x$ correspond to zero.
\end{theorem*}

\noindent\textit{Remark:} Consider two slices $S, S'$ in $M$ with, let us say, $\dimens S \leqq \dimens S'$. Suppose $g \in \G(S,S')$. Put $x = \s(g) \in S$ and ${x'} = \t(g) \in S'$. It is not difficult to see that there is a smooth target section $\tau: B \to \mca[1]\G$ defined over some open neighbourhood $B$ of ${x'}$ in $S'$ such that $\tau({x'}) = g$ and the composite map ${\s\circ\tau}$ induces a submersion of $B$ onto an open neighbourhood of $x$ in $S$. Thus, when \G\ is proper, it follows from the preceding theorem that for each point $x \in M$ there are a finite dimensional linear representation $G \to \GL(V)$ of a compact Lie group $G$ and a \G-invariant open neighbourhood $U$ of $x$ in $M$ for which there exists a Morita equivalence $\iota: {G \ltimes V} \into \G|_U$ such that $\mca[0]\iota: V \into U$ is an embedding of manifolds mapping the origin of $V$ to $x$.
\end{paragrafo}

\sezione{The Language of Tensor Categories}\label{ii}
This section consists of two parts. The first one contains a rather concise account of some basic standard categorical notions, a detailed exposition of which can be found for example in \cite{MacLane71,DeMi82,De90}. In the second part, and precisely from \S\ref{TL} onwards, we establish a couple of fundamental propositions for later use in the proof of our reconstruction theorem (Theorem \ref{xx.18}).

\bigskip

\noindent A \textit{tensor structure} on a category \Kt\ consists of the following data:
\begin{equazione}\label{TC1}
\text{a bifunctor~} \otimes: {\Kt\times\Kt} \longto \Kt\text, \quad \text{a distinguished object~} \TU \in \Ob(\Kt)
\end{equazione}
and a list of natural isomorphisms, called \textit{\mbox{ACU} constraints,}
\begin{equazione}\label{TC2}
\begin{split}
\begin{array}{c}
\alpha_{R,S,T}: {R\otimes ({S\otimes T})} \isoto {({R\otimes S})\otimes T}\text,
\\[\medskipamount]
\gamma_{R,S}: {R\otimes S} \isoto {S\otimes R}\text,
\\[\medskipamount]
\lambda_R: R \isoto {\TU\otimes R} \qquad \text{and} \qquad \rho_R: R \isoto {R\otimes\TU}
\end{array}
\end{split}
\end{equazione}
satisfying MacLane's ``{coherence conditions}'' (cf.\ for example MacLane \cite{MacLane71}, pp.~157~ff.\ and especially p.~180 for a detailed account). A \textit{tensor category} is a category endowed with a tensor structure. In the terminology of \cite{MacLane71}, the present notion corresponds to that of ``{symmetric monoidal category}''. The natural isomorphism $\alpha$, resp.\ $\gamma$ is called the associativity, resp.\ commutativity constraint; $\lambda$ and $\rho$ are called the (tensor) unit constraints.

In practice, we shall deal exclusively with ``complex'' tensor categories. Recall that a $k$-linear category, where $k$ is any number field, is a category \Kt\ whose hom-sets are each endowed with a structure of vector space over $k$ with respect to which composition of arrows is bilinear. One also says that \Kt\ is a category endowed with a $k$-linear structure. A {$k$-linear tensor category} is a tensor category endowed with a $k$-linear structure such that the bifunctor $\otimes$ is bilinear. From now on, in this paper, ``{tensor category}'' will mean ``{additive \nC-linear tensor category}''. Thus, in particular, there will be a zero object and for all objects $R,S$ there will be a direct sum ${R\oplus S}$.

Let \Kt[C], \Kt[C'] be tensor categories. A \textit{tensor functor} of \Kt[C] into \Kt[C'] is obtained by attaching, to an ordinary functor $F: \Kt[C] \to \Kt[C']$, two isomorphisms
\begin{equazione}\label{TFC}
\begin{split}
\begin{array}{l}
\tau_{R,S}: {F(R) \otimes F(S)} \isoto F({R\otimes S}) \quad \text{(natural in $R,S$)} \quad \text{and}
\\[\medskipamount]
\upsilon: \TU' \isoto F(\TU)\text,
\end{array}
\end{split}
\end{equazione}
called tensor functor constraints, which are required to satisfy certain conditions expressing their compatibility with the \textit{ACU} constraints of the tensor categories \Kt[C] and \Kt[C']. The reader is referred to \textit{loc.~cit.} for a discussion of these conditions. Recall that a functor of $k$-linear categories is said to be linear if the induced maps of hom-sets are $k$-linear. A linear functor between additive $k$-linear categories will preserve zero objects and direct sums. We agree that an assumption of linearity on the functor $F: \Kt[C] \to \Kt[C']$ is part of our definition of the notion of tensor functor.

Let $F,F'$ be tensor functors of \Kt[C] into \Kt[C']. A natural transformation $\lambda: F \to F'$ is said to be \textit{tensor preserving} if the following diagrams commute:
\begin{equazione}\label{TPNT}
\begin{split}
\xymatrix@C=25pt{{F(R)\otimes F(S)}\ar[d]^{\tau_{R,S}}\ar[rr]^-{\lambda(R)\otimes \lambda(S)} & & {F'(R)\otimes F'(S)}\ar[d]^{\smash{\tau'}_{R,S}} & \TU'\ar[d]^(.43){\upsilon\phantom'}\ar@{=}[r] & \TU'\ar[d]^(.43){\upsilon'} \\ F({R\otimes S})\ar[rr]^-{\lambda({R\otimes S})} & & F'({R\otimes S}) & F(\TU)\ar[r]^{\lambda(\TU)} & F'(\TU)\text.\!\!}
\end{split}
\end{equazione}
The collection of all tensor preserving natural transformations $F \to F'$ will be denoted by $\Hom^\otimes(F,F')$. Note that any natural transformation of $F$ into $F'$ is necessarily additive i.e.\ satisfies $\lambda({R\oplus S}) = {\lambda(R) \oplus \lambda(S)}$.

\begin{paragrafo}[Tensor* categories.]\label{T*C}
By an {anti-involution} on a tensor category \Kt\ we mean an anti-linear tensor functor
\begin{equazione}\label{T*C1}
*: \Kt \to \Kt\text, \quad R \mapsto R^*
\end{equazione}
with the property that there exists a tensor preserving natural isomorphism
\begin{equazione}\label{T*C2}
\iota_R: R^{**} \isoto R \quad \text{with} \quad \iota(R^*) = \iota(R)^*\text.
\end{equazione}
By fixing one such isomorphism once and for all, one obtains a mathematical structure which shall here be referred to as a \textit{tensor* category.} A morphism of tensor* categories, or \textit{tensor* functor,} is obtained by attaching, to an ordinary (linear) tensor functor $F$, a tensor preserving natural isomorphism
\begin{equazione}\label{T*F1}
\xi_R: F(R)^* \isoto F(R^*)
\end{equazione}
such that the following diagram commutes:
\begin{equazione}\label{T*F2}
\begin{split}
\xymatrix@R=13pt@C=15pt{F(R)^{**}\ar[r]^-{\can^*}\ar[dr]_\can & F(R^*)^*\ar[r]^-\can & F(R^{**})\ar[dl]^{F(\can)} \\ & F(R)\text.\!\! &}
\end{split}
\end{equazione}
A morphism of tensor* functors, better known as a self-conjugate tensor preserving natural transformation, is defined to be a tensor preserving natural transformation making the following diagram commutative:
\begin{equazione}\label{T*PNT}
\begin{split}
\xymatrix@C=35pt{F(R)^*\ar[d]^{\xi_R}\ar[r]^-{\lambda(R)^*} & F'(R)^*\ar[d]^{{\xi'}_R} \\ F(R^*)\ar[r]^{\lambda(R^*)} & F'(R^*)\text.\!\!}
\end{split}
\end{equazione}
We write $\Hom^{\otimes,*}(F,F')$ for reference to such transformations.

\textit{Example: the category of vector spaces.} If $V$ is a complex vector space, we let $V^*$ denote the space obtained by retaining the additive structure of $V$ but changing the scalar multiplication into ${zv^*} = ({\overline zv})^*$. The star here indicates that a vector of $V$ is to be regarded as one of $V^*$. Since any linear map $f: V \to W$ also maps $V^*$ linearly into $W^*$, we can also regard $f$ as a linear map $f^*: V^* \to W^*$. Moreover, the unique linear map of ${V^*\otimes W^*}$ into $({V\otimes W})^*$ sending ${v^*\otimes w^*} \mapsto ({v\otimes w})^*$ is an isomorphism, and complex conjugation sets up a linear bijection between $\nC$ and $\nC^*$. This turns vector spaces into a complex tensor category $\SpV_\nC$ with $V^{**} = V$.

\textit{Example: vector bundles.} By a direct generalization of the preceding construction one obtains the tensor* category \VB{M} of smooth complex vector bundles (of locally finite rank) over a smooth manifold $M$. We shall identify \VB{\pt}, where \pt\ denotes the one-point manifold, with the tensor* category $\SpV_\nC$ introduced above. Notice that the pullback of vector bundles along a smooth mapping $f: N \to M$ determines an obvious tensor* functor $f^*$ of \VB{M} into \VB{N}.

Let \Kt\ be a tensor* category. By a ``{real structure on an object $R \in \Ob(\Kt)$}'' we mean an isomorphism $\mu: R \isoto R^*$ in \Kt\ such that the composite ${\mu^*\cdot\mu}$ equals the identity on $R$ modulo the canonical identification $R^{**} \can R$ provided by \refequ{T*C2}. Let $\nR(\Kt)$ denote the category whose objects are the pairs $(R,\mu)$ consisting of an object $R \in \Ob(\Kt)$ together with a real structure $\mu$ on $R$ and whose morphisms $a: (R,\mu) \to (S,\nu)$ are the morphisms $a: R \to S$ in \Kt\ such that ${\nu\cdot a} = {a^*\cdot\mu}$. Note that $\nR(\Kt)$ is naturally an \nR-linear category; further, there is an obvious induced tensor structure on it, which turns it into an \nR-linear tensor category.

As an example of this construction, observe that one has an obvious equivalence of (real) tensor categories between $\SpV_\nR$ and $\nR(\SpV_\nC)$: in one direction, to any real vector space $V$ one can assign the pair $({\nC\otimes V},{z\otimes v} \mapsto {\overline z\otimes v})$; conversely, any real structure $\mu: U \isoto U^*$ on a complex vector space $U$ determines the real eigenspace $U^\mu \subset U$ of all $\mu$-invariant vectors. There is an analogous equivalence between \VB[\nR]{M} and $\nR\bigl(\VB{M}\bigr)$, for each smooth manifold $M$.

Notice that any tensor* functor $F: \Kt[C] \to \Kt[D]$ induces an obvious \nR-linear tensor functor $\nR(F): \nR(\Kt[C]) \to \nR(\Kt[D])$. For any tensor* functors $F, G: \Kt[C] \to \Kt[D]$, the map $\lambda \mapsto \tilde\lambda$ where $\tilde\lambda(R,\mu) \equiv \lambda(R)$ provides a bijection
\begin{equazione}
\Hom^{\otimes,*}(F,G) \isoto \Hom^\otimes\bigl(\nR(F),\nR(G)\bigr)
\end{equazione}
between the self-conjugate tensor preserving transformations $F \to G$ and the tensor preserving transformations $\nR(F) \to \nR(G)$. Indeed, by exploiting the additivity of the tensor* category \Kt[C], one can construct a functor $\Kt[C] \to \nR(\Kt[C])$ which plays the same role as the functor that assigns a complex vector space the underlying real vector space: one chooses, for each pair $R,S$ of objects of \Kt[C], a direct sum ${R\oplus S}$; then the obvious isomorphism ${R\oplus R^*} \iso ({R\oplus R^*})^*$ is a real structure on ${R\oplus R^*}$. Observe that the functor $\nR(\Kt[C]) \to \Kt[C]$, $(R,\mu) \mapsto R$ has an analogous interpretation. One therefore sees that the formalism of tensor* categories is essentially equivalent to that of real tensor categories.
\end{paragrafo}

\noindent The next results are original. They will be put to use only in the final section of this paper, in the proof of the reconstruction theorem.

\begin{paragrafo}[Terminology.]\label{TL}
Let \Kt\ be a tensor* category and $F: \Kt \to \SpV_\nC$ a tensor* functor with values into (finite dimensional) complex vector spaces. Let $H$ be a topological group, and suppose a homomorphism of monoids is given
\begin{equazione}\label{TL1}
\pi: H \longto \End^{\otimes,*}(F)\text.
\end{equazione}
We say that $\pi$ is \textit{continuous} if for every object $R \in \Ob(\Kt)$ the induced representation
\begin{equazione}\label{TL2}
\pi_R: H \longto \End(F(R))
\end{equazione}
defined by setting $\pi_R(h) = {\pi(h)}(R)$ is continuous.
\end{paragrafo}

\begin{proposizione}\label{O.prp5}
Let $\Kt,F$, $H$ and $\pi$ be as in \refcnt{TL}, with $\pi$ continuous, and suppose, in addition, that $H$ is a compact Lie group. Assume that the following condition is satisfied:
\begin{elenco}
\item[\textsl{(*)}]for each pair of objects $R, S \in \Ob(\Kt)$, and for each $H$-equivariant homomorphism $A: F(R) \to F(S)$, there exists some arrow $R \xto a S$ in \Kt\ with $A = F(a)$.
\end{elenco}

Then the homomorphism $\pi$ is surjective; in particular, the monoid $\End^{\otimes,*}(F)$ is a group.
\end{proposizione}
\begin{proof}
Put $K = \kernel\pi \subset H$. This is a closed normal subgroup, because it coincides with the intersection ${\bigcap \kernel{\pi_R}}$ over all objects $R$ of \Kt. On the quotient $G = H/K$ there is a unique (compact) Lie group structure such that the quotient homomorphism $H \to G$ becomes a Lie group homomorphism. Every $\pi_R$ can be indifferently thought of as a continuous representation of $H$ or a continuous representation of $G$, and every linear map $A: F(R) \to F(S)$ is a morphism of $G$-modules if and only if it is a morphism of $H$-modules. Being continuous, every $\pi_R$ is also smooth.

We claim there exists an object $R_0$ of \Kt\ such that the corresponding $\pi_{R_0}$ is faithful as a representation of $G$. Indeed, by the compactness of the Lie group $G$, we can find $R_1, \ldots, R_\ell \in \Ob(\Kt)$ with the property that
\begin{equazione}\label{O.esp83}
{\kernel{\pi_{R_1}} \cap \cdots \cap \kernel{\pi_{R_\ell}}} = \{e\}\text,
\end{equazione}
where $e$ denotes the unit of $G$; compare \cite{BtD85}, p.~136. Then, if we set $R_0 = {R_1 \oplus \cdots \oplus R_\ell}$, the representation $\pi_{R_0}$ will be faithful because of the existence of an obvious isomorphism of $G$-modules
\begin{equazione}\label{O.equ12}
F({R_1 \oplus \cdots \oplus R_\ell}) \iso {F(R_1) \oplus \cdots \oplus F(R_\ell)}\text.
\end{equazione}

Now, it follows that the $G$-module $F(R_0)$ is a ``generator'' for the tensor* category $\ReV_\nC(G)$ of all continuous, finite dimensional, complex $G$-modules. Indeed, every irreducible such $G$-module $V$ embeds as a submodule of some tensor power ${F(R_0)^{\otimes k} \otimes (F(R_0)^*)^{\otimes\ell}}$, see for instance \cite{BtD85}, p.~137. Since each $\pi(h)$ is, by assumption, self-conjugate and tensor preserving, this tensor power will be naturally isomorphic to $F\left({{R_0}^{\otimes k} \otimes ({R_0}^*)^{\otimes\ell}}\right)$ as a $G$-module and hence for each object $V$ of $\ReV_\nC(G)$ there will be some object $R$ of \Kt\ such that $V$ embeds into $F(R)$ as a submodule.

Next, consider an arbitrary natural transformation $\lambda \in \End(F)$. Let $R$ be an object of the category \Kt, and let $V \subset F(R)$ be a submodule. The choice of a complement to $V$ in $F(R)$ determines an endomorphism of modules $P_V: F(R) \to V \into F(R)$ which, by the assumption \textsl{(*)}, comes from some endomorphism of $R$ in \Kt. This implies that the linear operators $\lambda(R)$ and $P_V$ on the space $F(R)$ commute with one another and, consequently, that $\lambda(R)$ maps the subspace $V$ into itself. We will usually omit any reference to $R$ and write simply $\lambda_V$ for the linear map that $\lambda(R)$ induces on $V$ by restriction. Note finally that, given another submodule $W \subset F(S)$ and an equivariant map $B: V \to W$, one always has
\begin{equazione}\label{O.equ11}
{B \cdot \lambda_V} = {\lambda_W \cdot B}\text.
\end{equazione}
To prove this identity, one first extends $B$ to an equivariant map $F(R) \to F(S)$ and then invokes \textsl{(*)} as before.

Let \forgetfulfunctor{G} denote the tensor* functor $\ReV_\nC(G) \longto \SpV_\nC$ that assigns each $G$-module the underlying vector space. As our next step, we will define an isomorphism of complex algebras
\begin{equazione}\label{x.8}
\theta: \End(F) \isoto \End(\forgetfulfunctor G)
\end{equazione}
so that the following diagram commutes
\begin{equazione}\label{O.equ13}
\begin{split}
\xymatrix@C=35pt@R=20pt{H\ar[d]_-{\text{proj.}}\ar[r]^-\pi & \End(F)\ar[d]_-\simeq^-\theta \\ G\ar[r]^-{\envelopemorph G} & \End(\forgetfulfunctor G)\text,\!\!}
\end{split}
\end{equazione}
where $\envelopemorph{G}(g)$ is, for each $g \in G$, the natural transformation of \forgetfulfunctor{G} into itself that assigns left multiplication by $g$ on $V$ to each $G$-module $V$. For each $G$-module $V$ there exists an object $R$ of \Kt\ together with an embedding $V \into F(R)$, so we could define ${\theta(\lambda)}(V)$ as the restriction $\lambda_V$ of $\lambda(R)$ to $V$. Of course, it is necessary to check that this does not depend on the choices involved. Let two objects $R, S \in \Ob(\Kt)$ be given along with two equivariant embeddings of $V$ into $F(R),F(S)$. Since it is always possible to embed everything equivariantly into $F({R\oplus S})$ without affecting the induced $\lambda_V$, it is no loss of generality to assume $R=S$. Now, it follows from the remark \refequ{O.equ11} above that the two embeddings actually determine the same linear endomorphism of $V$. This shows that $\theta$ is well-defined. \refequ{O.equ11} also implies that $\theta(\lambda) \in \End(\forgetfulfunctor G)$. On the other hand put, for $\mu \in \End(\forgetfulfunctor G)$ and $R \in \Ob(\Kt)$, $\mu^F(R) = \mu(F(R))$. Then $\mu^F \in \End(F)$ and $\theta(\mu^F) = \mu$, because of the existence of embeddings of the form $V \into F(R)$ and because of the naturality of $\mu$. This shows that $\theta$ is surjective, and also injective since $\lambda(R) = \theta(\lambda)(F(R))$. Finally, it is straightforward to check that \refequ{O.equ13} commutes.

In order to conclude the proof it will be enough to check that $\theta$ induces a bijection between $\End^{\otimes,*}(F)$ and $\End^{\otimes,*}(\forgetfulfunctor G)$, for then our claim that $\pi$ is surjective will follow immediately from the commutativity of \refequ{O.equ13} and the classical \mbox{Tannaka}--\mbox{Kre\v\i n} duality theorem for compact groups (which says that \envelopemorph{G} establishes a bijection between $G$ and $\End^{\otimes,*}(\forgetfulfunctor G)$; see for example \cite{JoSt91} or \cite{BtD85} for a proof). This can safely be left to the reader.
\end{proof}

The argument we used in the foregoing proof to construct a ``generator'' tells us something interesting even in the noncompact case.
\begin{proposizione}\label{O.lem8}
Let \Kt\ and $F$ be as in \refcnt{TL}. Let $G$ be a Lie group, not necessarily compact, and let $\pi: G \to \End(F)$ be a faithful continuous homomorphism. Then there exists an object $R_0 \in \Ob(\Kt)$ such that $\kernel{\pi_{R_0}}$ is a discrete subgroup of $G$ or, equivalently, such that the representation
\begin{equazione}\label{x.14}
\pi_{R_0}: G \to \GL(F(R_0))
\end{equazione}
is faithful in some open neighbourhood of the unit of $G$.
\end{proposizione}
\begin{proof}
By induction.
\end{proof}

\sezione{Smooth Tensor Stacks}\label{iii}
In this section we shall introduce the language of smooth stacks of tensor* categories or, in short, {smooth tensor stacks.} We propose this language as a new foundational framework for the representation theory of groupoids. Nowadays, many concise accounts of the general theory of stacks are available; our own exposition is based on \cite{DeMi82} and \cite{Moe0212}. A fairly extensive treatment of the algebraic geometric theory can be found in \cite{LauM-B}.

\bigskip

\noindent\textit{Overall Conventions:} The capital letters $X,Y,Z$ denote $\C^\infty$ manifolds and the corresponding lower-case letters $x, x', \ldots, y$ etc.\ denote points on these manifolds. `\smooth[X]' stands for the sheaf of smooth functions on $X$; we sometimes omit the subscript. We refer to sheaves of \smooth[X]-modules also as sheaves of modules over $X$. For practical purposes, we need to consider manifolds which are possibly neither Hausdorff nor paracompact.

\begin{paragrafo}[Fibred tensor categories.]\label{Fib-T*C}
Fibred tensor categories shall be denoted by capital Gothic type variables. A fibred tensor category \stack C assigns, to each smooth manifold $X$, a tensor* category
\begin{equazione}\label{Fib-T*C1}
\stack[X]C = \bigl(\stack[X]C,\otimes_X,\TU_X,*_X\bigr)
\end{equazione}
or $\bigl(\stack[X]C,\otimes,\TU,*\bigr)$ for short\inciso{we omit subscripts when they are retrievable from the context}and, to each smooth map $X \xto f Y$, a tensor* functor
\begin{equazione}\label{Fib-T*C2}
f^*: \stack[Y]C \longto \stack[X]C
\end{equazione}
which we call pull-back along $f$. Moreover, for each pair of composable smooth maps $X \xto f Y \xto g Z$ and for each manifold $X$, any fibred tensor category provides self-conjugate tensor preserving natural isomorphisms
\begin{equazione}\label{N.i3}
\left\{\begin{aligned}
\delta&: {f^*\circ g^*} \isoto ({g\circ f})^*
\\
\varepsilon&: \Id_{\stack[X]C} \isoto {\id_X}^*
\end{aligned}\right.
\end{equazione}
which altogether are required to make the following diagrams commute:
\begin{equazione}\label{N.i4}
\begin{split}
\xymatrix@C=30pt@R=20pt{f^*g^*h^*\ar[d]^(.47){\delta\cdot h^*}\ar[r]^-{f^*\delta} & f^*(hg)^*\ar[d]^(.47)\delta & & {\id_X}^*f^*\ar[d]^(.47)\delta & f^*\ar[d]^(.47){f^*\varepsilon}\ar@{=}[dl]\ar[l]_(.34){\varepsilon\cdot f^*} \\ (gf)^*h^*\ar[r]^-\delta & (hgf)^* & & f^* & f^*{\id_Y}^*\text{~.}\mspace{-10.0mu}\ar[l]_(.5)\delta}
\end{split}
\end{equazione}
These are the only primitive data we need to introduce in order to be able to speak about smooth tensor stacks and representations of Lie groupoids. The latter concepts can\inciso{and will}be defined {\em canonically,} i.e.\ purely in terms of the fibred tensor category structure.
\end{paragrafo}

\begin{paragrafo}[Tensor prestacks.]\label{iii.2}
Let \stack C be an arbitrary fibred tensor category. Let $i_U: U \into X$ denote the inclusion of an open subset. We shall put, for every object $E$ and morphism $a$ of the category \stack[X]C, $E|_U = {i_U}^*E$ and $a|_U = {i_U}^*a$. More generally, we shall make use of the same abbreviations for the inclusion $i_S: S \into X$ of an arbitrary submanifold.

For each pair of objects $E, F \in {\Ob\,\stack[X]C}$, let \sheafhom[C]{X}{E}{F} denote the presheaf of complex vector spaces over $X$ defined as
\begin{equazione}\label{N.i5}
U \mapsto \Hom_{\stack[U]C}(E|_U,F|_U)
\end{equazione}
where the restriction map corresponding to an open inclusion $j: V \into U$ is given [obviously, up to canonical isomorphism] by $a \mapsto {j^*a}$. Now, the requirement that \stack C is a prestack means exactly that every such presheaf is in fact a sheaf. This entails, in particular, that one can glue any family of compatible local morphisms over $X$. One special case, namely the sheaf $\sections E = \sheafhom[C]{X}{\TU}{E}$, to which we shall refer as the \textit{sheaf of sections of $E$}, will be of major interest for us. For any open subset $U$, the elements of the vector space $\sections E(U)$ shall be referred to as \textit{sections of $E$ over $U$}.

Since a morphism $a: E \to F$ in \stack[X]{C} yields a morphism $\sections a: \sections E \to \sections F$ of sheaves of complex vector spaces over $X$ (by composing $\TU|_U \to E|_U \xto{a|_U} F|_U$), we obtain a canonical functor
\begin{equazione}\label{N.i6}
\sections{} = \sections[X]{}: \stack[X]{C} \longto \SheavesOfModules{\sheafconst[\mathnormal X]\nC}\text,
\end{equazione}
where \sheafconst[X]{\nC} denotes the constant sheaf over $X$ of value \nC.

This functor is certainly linear. Moreover, there is a canonical way to make it a ``pseudo-''tensor functor of the tensor category $\bigl(\stack[X]C,\otimes_X,\TU_X\bigr)$ into the category of sheaves of \sheafconst[X]{\nC}-modules (with the familiar tensor structure): a natural transformation $\tau_{E,F}: {\sections[X]E \otimes_{\sheafconst[X]\nC} \sections[X]F} \to \sections[X]{({E\otimes F})}$ arises, in the obvious manner, from the local pairings
\begin{equazione+}\label{N.i7}
{\sections E(U) \times \sections F(U)} \to \sections{({E\otimes F})}(U)\text, \quad (a,b) \mapsto {a\otimes b} & \text{[mod $\can$]}
\end{equazione+}
(which are bilinear with respect to locally constant coefficients), and a unit constraint $\upsilon: \sheafconst[X]\nC \to \sections[X]\TU$ may be defined as follows:
\begin{equazione}\label{N.i8}
{\scriptscriptstyle \left\{ \begin{array}{c} \mathrm{locally\:constant\:complex} \\ \mathrm{valued\:functions\:on}\: U \end{array} \right\}} \longto \sections{\TU}(U)\text, \quad z \mapsto {z \cdot \id|_U}\text;
\end{equazione}
it is easy to check that these morphisms of sheaves make the same diagrams which characterize tensor functor constraints commute.

One also has a natural morphism of sheaves of modules over $X$
\begin{equazione}\label{xi35}
(\sections[X]E)^* \longto \sections[X]{(E^*)}
\end{equazione}
defined by means of the anti-involution and the obvious related canonical isomorphisms. Since $\zeta^{**} = \zeta$ [up to canonical isomorphism], it follows at once that \refequ{xi35} is a natural {\em isomorphism} for any tensor prestack. In fact, \refequ{xi35} is an isomorphism of pseudo-tensor functors viz.\ it is compatible with the natural transformations \refequ{N.i7} and \refequ{N.i8}.
\end{paragrafo}

\begin{paragrafo}[Fibres of an object.]\label{iii.3}
Note that for $X = \pt$ (where \pt\ is the one-point manifold) one may regard \refequ{N.i6} as a {pseudo-tensor* functor} of \stack[\pt]C into complex vector spaces. We put, for all objects $E \in {\Ob\,\stack[\pt]C}$,
\begin{equazione}\label{N.i10}
E_* = (\sections[\pt]E)(\pt)
\end{equazione}
(so this is a complex vector space) and do the same for morphisms. Now, as a part of the forthcoming definition of the general notion of smooth tensor prestack, we impose the following requirement: \textsl{the morphism of sheaves \refequ{N.i8} is an isomorphism for $X = \pt$.} Then there is a canonical isomorphism
\begin{equazione}\label{N.i11}
\nC \can \TU_*
\end{equazione}
of complex vector spaces.

For any object $E \in {\Ob\,\stack[X]C}$, we define the \textit{fibre of $E$ at $x$} to be the complex vector space $E_x = (x^*E)_*$. We use the same name for the point $x$ and for the (smooth) mapping $\pt \to X, \pt \mapsto x$, so that $x^*$ is just our ordinary notation \refequ{Fib-T*C2} for the pull-back, $x^*E$ belongs to \stack[\pt]C and we can make use of the notation \refequ{N.i10}. Similarly, whenever $a: E \to F$ is a morphism in \stack[X]{C}, we let $a_x: E_x \to F_x$ denote the linear map $(x^*a)_*$. Since $\text- \mapsto (\text-)_x$ is by construction the composite of two pseudo-tensor* functors, it itself may be regarded as a pseudo-tensor* functor. If in particular we apply this to a local section $\zeta \in \sections E(U)$ and make use of the canonical identification \refequ{N.i11}, we get, for $u \in U$, a linear map
\begin{equazione}\label{N.i12}
\nC \can \TU_* \can ({u^*\,\TU|_U})_* \xto{(u^*\zeta)_*} ({u^*\,E|_U})_* \can (u^*E)_* = E_u
\end{equazione}
which corresponds to a vector $\zeta(u) \in E_u$ to be called the \textit{value of $\zeta$ at $u$}. One has the intuitive formula
\begin{equazione}\label{N.i13}
a_u(\zeta(u)) = \bigl[(\sections a)(U)(\zeta)\bigr](u)\text.
\end{equazione}
Notice finally that the vectors ${\zeta(u) \otimes \eta(u)}$ and $({\zeta\otimes\eta})(u)$ correspond to one another via the canonical map ${E_u \otimes_\nC F_u} \to ({E\otimes_X F})_u$.
\end{paragrafo}

\begin{paragrafo}[Smooth tensor prestacks.]\label{iii.4}
Let $\TU_X$ denote the tensor unit of \stack[X]{C}, and let $x$ be a point of the manifold $X$. Under the assumption \refequ{N.i11}, one can use the composite linear isomorphism $\nC \can (\TU_\pt)_* \can ({x^*\TU_X})_* \equiv (\TU_X)_x$ to define a canonical homomorphism of complex algebras
\begin{equazione}\label{N.i14}
\End_{\stack[X]C}(\TU_X) \longto \left\{ \mathrm{complex\:functions\:on\:} X \right\}\text, \quad e \mapsto \tilde e
\end{equazione}
by putting $\tilde e(x)$ = the complex number such that the linear map \gm{scalar multiplication by $\tilde e(x)$} (of \nC\ into itself) corresponds to $e_x: (\TU_X)_x \to (\TU_X)_x$. We shall say that a tensor prestack \stack C is \textit{smooth} if it verifies \refequ{N.i11} and if the homomorphism \refequ{N.i14} determines a bijective correspondence between $\End_{\stack[X]C}(\TU_X)$ and the subalgebra of all smooth functions on $X$:
\begin{equazione}\label{N.i15}
\End(\TU_X) \can \C^\infty(X)\text.
\end{equazione}

When a tensor prestack \stack C is smooth, it is possible to endow each space $\Hom_{\stack[X]C}(E,F)$ with a canonical $\C^\infty(X)$-module structure compatible with the multiplication by locally constant functions, since $\Hom(E,F)$ has an obvious $\End(\TU_X)$-module structure. Accordingly, $\sheafhom[C]{X}{E}{F}(U)$ inherits a canonical structure of $\C^\infty(U)$-module for every open subset $U \subset X$, which turns \sheafhom[C]{X}{E}{F} into a sheaf of \smooth[X]-modules. This is true, in particular, of the sheaf of smooth sections \sections[X]{E}. It is also readily seen that each morphism $a: E \to F$ induces a morphism $\sections[X]a: \sections[X]E \to \sections[X]F$ of sheaves of \smooth[X]-modules. Thus, one gets a $\C^\infty(X)$-linear functor
\begin{equazione}\label{N.i17}
\sections[X]{}: \stack[X]{C} \longto \SheavesOfModules{\smooth[\mathnormal X]}\text.
\end{equazione}

The operations \refequ{N.i7} and \refequ{N.i8} may now be used to define morphisms of sheaves of \smooth[X]\nobreakdash-modules
\begin{equazione}\label{N.i18}
\left\{\begin{aligned}
\tau&: {\sections[X]E \otimes_{\smooth[X]} \sections[X]F} \to \sections[X]{({E\otimes F})}
\\
\upsilon&: \smooth[X] \to \sections[X]\TU\text.
\end{aligned}\right.
\end{equazione}
The morphism $\tau = \tau_{E,F}$ is natural in the variables $E,F$ and, along with $\upsilon$, turns \refequ{N.i17} into a pseudo-tensor functor of \stack[X]C into the category of sheaves of \smooth[X]-modules. This is closer than \refequ{N.i6} to being a genuine tensor functor, in that the morphism $\upsilon$ is now an isomorphism of sheaves of \smooth[X]-modules.

Consider next a smooth mapping of manifolds $f: X \to Y$. Suppose that $U \subset X$ and $V \subset Y$ are open subsets with $f(U) \subset V$, and let $f_U$ denote the induced mapping of $U$ into $V$. For any object $F$ of the category \stack[Y]C, a correspondence of local smooth sections
\begin{equazione}\label{N.i19}
(\sections[Y]F)(V) \longto \sections[X]{({f^*F})}(U)\text, \quad \eta \mapsto {\eta\circ f}
\end{equazione}
can be obtained by defining ${\eta\circ f}$ as the following composite arrow:
\begin{equazione}\label{N.i20}
\TU_X|_U \can ({f^*\TU_Y})|_U \can {f_U}^*(\TU_Y|_V) \xto{{f_U}^*(\eta)} {f_U}^*(F|_V) \can ({f^*F})|_U\text.
\end{equazione}
For $U$ fixed, and $V$ variable, the maps \refequ{N.i19} constitute an inductive system indexed over the inclusions $V \supset V' \supset f(U)$ and hence they yield, on passing to the limit, a morphism of sheaves of \smooth[X]-modules
\begin{equazione}\label{N.i21}
f^*(\sections[Y]F) \longto \sections[X]{({f^*F})}\text,
\end{equazione}
where $f^*(\sections[Y]F)$ is the ordinary pullback of sheaves of modules over a smooth manifold. The morphism \refequ{N.i21} is natural in $F$. It is also a morphism of pseudo-tensor functors, in other words it is tensor preserving. Notice that the vectors $\eta(f(x)) \in F_{f(x)}$ and $({\eta\circ f})(x) \in (f^*F)_x$ correspond to one another via the canonical isomorphism of vector spaces
\begin{equazione}\label{N.i22}
({f^*F})_x = ({x^*f^*F})_* \can ({f(x)^*F})_* = F_{f(x)}\text.
\end{equazione}
\end{paragrafo}

\begin{paragrafo}[Flat maps.]\label{iii.5}
It will be convenient to regard open coverings of a manifold as smooth maps. We say that a smooth map $p: X' \to X$ is \textit{flat} if it is surjective and it restricts to an open embedding $p_{U'}: U' \into X$ on each connected component $U'$ of $X'$. We may think of $p$ as codifying the open covering of $X$ given by the connected components of $X'$. A {refinement} of $X' \xto p X$ is obtained by composing $p$ backwards with another flat mapping $X'' \xto{p'} X'$. If $p$ is flat then for any smooth map $f: Y \to X$ the pullback
\begin{equazione}\label{N.i24}
{Y \times_X X'} = \{(y,x'): f(y)=p(x')\}
\end{equazione}
exists in the category of $\C^\infty$ manifolds and $\pr_1: {Y \times_X X'} \to Y$ is also flat. When $f$ itself is flat this construction yields a common refinement. For any flat map $p: X' \to X$, put
\begin{equazione}\label{N.i25}
X'' \bydef {X' \times_X X'} = \{(x'_1,x'_2): p(x'_1)=p(x'_2)\}
\end{equazione}
and let $p_1, p_2: X'' \to X'$ be the projections. Also put
\begin{equazione}\label{N.i26}
X''' \bydef {X'\times_X X'\times_X X'} = \{(x'_1,x'_2,x'_3): p(x'_1)=p(x'_2)=p(x'_3)\}
\end{equazione}
and let $p_{12}: X''' \to X''$ be the map $(x'_1,x'_2,x'_3) \mapsto (x'_1,x'_2)$; the maps $p_{23}$ and $p_{13}$ shall be defined likewise.
\end{paragrafo}

\begin{paragrafo}[Smooth tensor stacks.]\label{iii.6}
A descent datum for a smooth tensor prestack \stack C relative to a flat mapping $p: X' \to X$ is a pair $(E',\theta)$ consisting of an object $E'$ of the category \stack[X']C and an isomorphism $\theta: {{p_1}^*E'} \isoto {{p_2}^*E'}$ in the category \stack[X'']{C} such that ${p_{13}}^*(\theta) = {{p_{12}}^*(\theta) \circ {p_{23}}^*(\theta)}$ [mod $\can$]. A morphism $a': (E',\theta) \to (F',\xi)$ of descent data relative to $p$ is an arrow $a': E' \to F'$ compatible, in the obvious sense, with the isomorphisms $\theta$ and $\xi$. Descent data for \stack C relative to $p$ and their morphisms form a category \Des[C]{X'/X}. There is an obvious functor
\begin{equazione}\label{N.i31}
\stack[X]C \longto \Des[C]{X'/X}\text, \quad E \mapsto ({p^*E},\phi_E)
\end{equazione}
where $\phi_E$ is the canonical isomorphism ${p_1}^*({p^*E}) \can {({p\circ p_1})^*E} = {({p\circ p_2})^*E} \can {p_2}^*({p^*E})$. If this canonical functor is an equivalence of categories for every flat mapping $p: X' \to X$, one says that \stack C is a stack.

For our purposes, the condition that the functor \refequ{N.i31} be an equivalence of categories for every flat map $X' \to X$ can be relaxed to some extent. In fact, it is sufficient to require it of all flat maps $X' \to X$ over a Hausdorff, paracompact $X$. We propose to use the term `parastack' for the weaker notion. We will often be sloppy and use `stack' as a synonym to `parastack'.
\end{paragrafo}

\begin{paragrafo}[Locally trivial objects.]\label{iii.7}
Let \stack C be a smooth tensor prestack. An object $E \in {\Ob\,\stack[X]C}$ will be called trivial if there exists some $V \in {\Ob\,\stack[\pt]C}$ for which one can find an isomorphism $E \iso {{c_X}^*V}$ in \stack[X]{C} where $c_X: X \to \pt$ denotes the collapse map. Any such pair $(V,\iso)$ will be said to constitute a trivialization of $E$.

For an arbitrary manifold $X$, let \V[C]X denote the full subcategory of \stack[X]C formed by the locally trivial objects of locally finite rank: $E \in {\Ob\,\stack[X]C}$ is an object of \V[C]X if and only if one can cover $X$ with open subsets $U$ such that there exists in each \stack[U]{C} a trivialization $E|_U \iso {\TU_U \oplus \cdots \oplus \TU_U}$. The operation $X \mapsto \V[C]X$ determines a fibred tensor subcategory of \stack C. In fact, one gets a smooth tensor prestack $\mathit V^{\stack C}$ which is a parastack resp.\ a stack if so is \stack C.

The tensor* category \V[C]X closely relates to that of smooth complex vector bundles over $X$. Clearly, every object $E \in \V[C]X$ yields a vector bundle $\tilde E = \{(x,e): x \in X, e \in E_x\}$ over $X$. The operation $E \mapsto \tilde E$ defines a faithful tensor* functor of \V[C]{X} into \VB{X} which is an equivalence of tensor* categories when \stack C is a parastack and $X$ is paracompact or when \stack C is a stack.
\end{paragrafo}

\sezione{Foundations of Representation Theory}\label{iv}
In this section, we develop the representation theory of Lie groupoids within the framework described in \S\ref{iii}. A peculiarity of the notion of representation we shall be considering here is its dependence on a `type': the construction of our theory necessitates the preliminary choice of an arbitrary smooth tensor stack \stack T (the type). We shall assume that such a choice has been made and we shall regard \stack T as fixed throughout the section. The definitions below are directly inspired by \cite{De90}.

Let \G\ be a Lie groupoid over a manifold $M$. We start by constructing the category of representations ``of type \stack T'' of \G. Define the objects of \R[T]{\G}, or briefly, \R{\G}, to be the pairs $(E,\varrho)$ consisting of an object $E$ of \stack[M]{T} and an arrow $\varrho$ in \stack[\G]{T}
$$%
\s^*E \xto{\:\:\varrho\:\:} \t^*E\text,
$$%
where $\s, \t: \G \to M$ denote the source resp.\ target map of \G, such that the appropriate conditions for $\varrho$ to be an action, in other words for $\varrho$ to be compatible with the groupoid structure, are satisfied (modulo the appropriate canonical isomorphisms):
\begin{elenco}
\item ${\u^* \varrho} = \id_E$, where $\u: M \to \G$ denotes the unit section of \G;
\item ${\c^* \varrho} = {{{\p_0}^* \varrho} \cdot {{\p_1}^* \varrho}}$, where $\mca[2]\G = \mca\G$ denotes the manifold of composable arrows, $\c: \mca[2]\G \to \G$ the composition map and $\p_0, \p_1: \mca[2]\G \to \G$ the left resp.\ right projection.
\end{elenco}
Observe that the conditions i) and ii) imply that the arrow $\s^*E \xto\varrho \t^*E$ is invertible in the category \stack[\G]{T}. We shall refer to the objects of \R{\G} also as \textit{\G-actions} (\textit{of type \stack T}). Define the morphisms of \G-actions $(E,\varrho) \to (E',\varrho')$ to be the arrows $a: E \to E'$ in \stack[M]{T} which fulfill the condition
\begin{equazione}\label{N.ii4}
{\t^*a \cdot \varrho} = {\varrho' \cdot \s^*a}\text.
\end{equazione}
The category \R{\G}, endowed with the \nC-linear structure inherited from \stack[M]{T}, is clearly additive.

\begin{paragrafo}[Tensor* structure.]\label{iv.1}
For any \G-actions $R = (E,\varrho)$ and $S = (F,\sigma)$ we put ${R\otimes S} = ({E\otimes F},{\varrho\otimes\sigma})$ where ${\varrho\otimes\sigma}$ stands for
\begin{equazione}\label{N.ii10}
\s^*({E\otimes F}) \can {{\s^*E} \otimes {\s^*F}} \xto{\:\:\varrho\,\otimes\,\sigma\:\:} {{\t^*E} \otimes {\t^*F}} \can \t^*({E\otimes F})\text.
\end{equazione}
It is easy to recognize that the pair ${R\otimes S}$ itself is a \G-action. Further, if $(E,\varrho) \xto a (E',\varrho')$ and $(F,\sigma) \xto b (F',\sigma')$ are morphisms of \G-actions then so will be ${a\otimes b}: {R\otimes S} \to {R'\otimes S'}$. We define the tensor unit of \R{\G} as the pair $(\TU_M,\can)$, where $\TU_M$ denotes the tensor unit of \stack[M]{T} and $\can$ stands for
\begin{equazione}\label{N.ii11}
{\s^*\TU_M} \can \TU_\G \can {\t^*\TU_M}\text.
\end{equazione}
The \textit{ACU} constraints for the tensor product on the base category \stack[M]{T} provide analogous constraints for the tensor product on \R{\G}. There is of course also an inherited canonical anti-involution.

The forgetful functor
\begin{equazione}\label{N.ii5}
\forgetfulfunctor[T]{\G}: \R[T]{\G} \longto \stack[M]{T}\text, \quad (E,\varrho) \mapsto E
\end{equazione}
(or \forgetfulfunctor{\G}, for short) is \nC-linear and faithful. By construction, it is a strict tensor* functor of \R{\G} into \stack[M]{T}. We shall refer to this functor as the \textit{standard forgetful functor} (\textit{of type \stack T}) associated with \G.
\end{paragrafo}

\begin{paragrafo}[Pullback along a homomorphism.]\label{iv.2}
Let $\varphi: \G \to \H$ be a homomorphism of Lie groupoids and let $M \xto f N$ be the smooth map induced by $\varphi$ on the base manifolds.

Suppose $(F,\sigma) \in \R[T]\H$. Consider the morphism\inciso{which we still denote by ${\varphi^*\sigma}$, allowing some notational abuse}defined as follows:
\begin{equazione}\label{N.ii15}
{\s_{\G}}^*({f^*F}) \can {\varphi^*{\s_{\H}}^*F} \xto{\:\varphi^*\sigma\:} {\varphi^*{\t_{\H}}^*F} \can {\t_{\G}}^*({f^*F})\text.
\end{equazione}
The identities ${f\circ\s_{\G}} = {\s_{\H}\circ\varphi}$ etc.\ account, of course, for the existence of the canonical isomorphisms occurring in \refequ{N.ii15}. It is routine to check that the pair $({f^*F},{\varphi^*\sigma})$ constitutes an object of the category \R[T]\G\ and that if $(F,\sigma) \xto b (F',\sigma')$ is a morphism of \H\nobreakdash-actions then ${f^*b}$ is a morphism of $({f^*F},{\varphi^*\sigma})$ into $({f^*F'},{\varphi^*\sigma'})$ in \R[T]\G. Hence we get a functor
\begin{equazione}\label{N.ii16}
\varphi^*: \R[T]\H \longto \R[T]\G\text,
\end{equazione}
which we agree to call the \textit{inverse image} (or \textit{pull-back}) \textit{along $\varphi$}.

The constraints associated with the tensor functor $f^*$
\begin{equazione}\label{N.ii17}
\left\{\begin{aligned}
& \TU_M \isoto {f^*\TU_N}
\\
& {{f^*F} \otimes {f^*F'}} \isoto f^*({F \otimes F'})
\end{aligned}\right.
\end{equazione}
function as isomorphisms of \G-actions $\TU \isoto \varphi^*(\TU)$ and ${\varphi^*(S) \otimes \varphi^*(S')} \isoto \varphi^*({S\otimes S'})$ for all $S, S' \in \R\H$ with $S = (F,\sigma)$ and $S' = (F',\sigma')$. A fortiori, these isomorphisms are natural and provide appropriate tensor functor constraints for $\varphi^*$, thus making $\varphi^*$ a tensor functor of the tensor category \R\H\ into the tensor category \R\G.

Let $\G \xto\varphi \H \xto\psi \K$ be two composable homomorphisms of Lie groupoids and let $X \xto{\varphi_0} Y \xto{\psi_0} Z$ denote the respective base maps. Note that, for an arbitrary action $T = (G,\tau) \in \R{\K}$, the canonical isomorphism ${{\varphi_0}^* \smash{\psi_0}^* G} \can {({\psi_0 \circ \varphi_0})^* G} = {\smash{({\psi \circ \varphi})_0}^* G}$ is an isomorphism between $\varphi^*({\psi^* T})$ and ${({\psi \circ \varphi})^* T}$ in the category \R{\G}. Hence we get an isomorphism of tensor functors
\begin{equazione}\label{xiv.4}
{\varphi^* \circ \psi^*} \xto{\:\simeq\:} ({\psi \circ \varphi})^*\text.
\end{equazione}
\end{paragrafo}

\begin{paragrafo}[Natural transformations.]\label{iv.3}
Recall that a transformation $\tau: \varphi_0 \isoto \varphi_1$ between two Lie groupoid homomorphisms $\varphi_0, \varphi_1: \G \to \H$ is a smooth mapping $\tau$ of the base manifold $M$ of \G\ into the manifold of arrows of \H\ such that $f_0(x) \xto{\tau(x)} f_1(x)$ for all $x \in M$ and
\begin{equazione}\label{N.ii19}
{\varphi_1(g) \cdot \tau(x)} = {\tau(x') \cdot \varphi_0(g)}
\end{equazione}
for all $g \in \mca[1]{\G}, g: x \to x'$. Suppose an action $S = (F,\sigma) \in \R{\H}$ is given. Then we can apply $\tau^*$ to the isomorphism ${\s^*F} \xto\sigma {\t^*F}$ to obtain an isomorphism ${f_0^*F} \to {f_1^*F}$ in the category \stack[M]{T}
\begin{equazione}\label{N.ii20}
{f_0^*F} \can {\tau^*\s^*F} \xto{\:\tau^*\sigma\:} {\tau^*\t^*F} \can {f_1^*F}\text.
\end{equazione}
By expressing \refequ{N.ii19} as an identity between suitable smooth maps, one sees that \refequ{N.ii20} is actually an isomorphism of \G-actions (between ${\varphi_0^*S}$ and ${\varphi_1^*S}$). Thus, we obtain an isomorphism of tensor functors $\varphi_0^* \simeq \varphi_1^*$.
\end{paragrafo}

\begin{paragrafo}[Morita equivalences.]\label{iv.4}
We observe next that the inverse image functor $\varphi^*: \R\H \to \R\G$ associated with a Morita equivalence $\varphi: \G \to \H$ is an equivalence of tensor categories.\footnote{Recall that a tensor functor $\Phi: \Kt[C] \to \Kt[D]$ is said to be a tensor equivalence in case there exists a tensor functor $\Psi: \Kt[D] \to \Kt[C]$ for which there are tensor preserving natural isomorphisms ${\Psi\circ\Phi} \simeq \Id_{\Kt[C]}$ and ${\Phi\circ\Psi} \simeq \Id_{\Kt[D]}$.} Clearly, this is tantamount to saying that $\varphi^*$ is a categorical equivalence. Although the procedure to obtain a quasi-inverse $\varphi_!$ follows a well-known pattern, we review it for the reader's convenience. In fact, we know of no adequate standard reference for this precise argument.

The condition that the map \refequ{ME2} be a surjective submersion will of course be satisfied when \mca[0]{\varphi} itself is a surjective submersion. As a first step, we show how the task of constructing a quasi-inverse may be reduced to the special case where \mca[0]{\varphi} is precisely a surjective submersion. To this end, consider the weak pullback (see \cite{MoeMrc03}, pp.~123--132)
\begin{equazione}\label{xiv.10}
\begin{split}
\xymatrix@C=40pt{\mathcal P\ar[d]^\chi\ar[r]^\psi & \G\ar[d]^\varphi_{}="b" \\ \H\ar@{=}[r]^(.4)\quad="a"\ar@{=>}@/^.5pc/"a";"b"^\tau & \H\text.\!\!}
\end{split}
\end{equazione}
Let $P$ be the base manifold of the Lie groupoid $\mathcal P$. It is well-known (\textit{ibid.}\ p.~130) that the Lie groupoid homomorphisms $\psi$ and $\chi$ are Morita equivalences with the property that the respective base maps $\mca[0]\psi: P \to M$ and $\mca[0]\chi: P \to N$ are surjective submersions. Now, if we prove that $\psi^*$ and $\chi^*$ are categorical equivalences then, since by (\ref{xiv.4}) and the remarks contained in \S\ref{iv.3} we have natural isomorphisms
\begin{equazione}\label{xiv.11}
\chi^* \xto{\:\simeq\:} ({\varphi\circ\psi})^* \xfrom{\:\simeq\:} {\psi^* \circ \varphi^*}\text,
\end{equazione}
the same will be true of $\varphi^*$.

From now on, we work under the hypothesis that the Morita equivalence $\varphi$ determines a surjective submersion $f: M \to N$ on the base manifolds. This being the case, there exists an open cover of the manifold $N = \txtcup{i\in I}{}{V_i}$ by open subsets $V_i$ such that for each of them one can find a smooth section $s_i: V_i \into M$ to $f$. We fix such a cover and such sections once and for all.

Let an arbitrary object $R = (E,\varrho) \in \R{\G}$ be given. For each $i \in I$ one can take the pullback $E_i \equiv {{s_i}^*E} \in \stack[V_i]{T}$. Fix a couple of indices $i, j \in I$. Then, since \refequ{ME1} is a pullback diagram, for each $y \in {V_i\cap V_j}$ there is exactly one arrow $g(y): s_i(y) \to s_j(y)$ such that $\varphi(g(y)) = y$. More precisely, let $y \mapsto g(y) = g_{ij}(y)$ be the smooth mapping defined as the unique solution to the following universal problem (in the $\C^\infty$ category)
\begin{equazione}\label{xiv.12}
\begin{split}
\xymatrix@C=30pt{V_{ij}\ar@/_1.5pc/[ddr]_(.4){(s_i,s_j)}\ar@{-->}[dr]^(.55){g_{ij}}\ar@/^1.3pc/[drr]^(.6){\u|V_{ij}} \\ & \G\ar[d]^{(\s,\t)}\ar[r]^\varphi & \H\ar[d]^{(\s,\t)} \\ & {M\times M}\ar[r]^{f\times f} & {N\times N}\text,\!\!}
\end{split}
\end{equazione}
where $\u: N \to \H$ denotes the unit section and $V_{ij} = {V_i\cap V_j}$. Then, putting $E_i|_j = E_i|_{V_i\cap V_j}$ and $E_j|_i = E_j|_{V_i\cap V_j}$, one may pull the action $\varrho$ back along the map $g_{ij}$ so as to get an isomorphism $\theta_{ij}: E_i|_j \isoto E_j|_i$
\begin{equazione+}\label{xiv.13'}
\theta_{ij} = {{g_{ij}}^*\varrho} & \text{[mod $\can$]}
\end{equazione+}
in the category \stack[V_{ij}]{T}. Next, from the obvious remark that for an arbitrary third index $k \in I$ one has $g_{ik}|_j = {\c\circ (g_{jk}|_i,g_{ij}|_k)}$, where $g_{ik}|_j$ denotes the restriction of $g_{ik}$ to $V_{ijk}$, and from the multiplicative axiom ii) for $\varrho$, it follows that the system of isomorphisms $\{\theta_{ij}\}$ constitutes a ``cocycle'' or ``descent datum'' for the family $\{E_i\}_{i\in I} \in \xstack[{\smash[b]{\, \txtcoprod{i\in I}{}{V_i} \,}}]{T}$ relative to the flat mapping $\txtcoprod{i\in I}{}{V_i} \to N$. Since $N$ is a paracompact manifold and \stack T is a smooth parastack, there exist an object ${\varphi_!E}$ of \stack[N]{T} and a system of isomorphisms $\theta_i: ({\varphi_!E})|_i \equiv ({\varphi_!E})|_{V_i} \isoto E_i$ in \stack[V_i]{T} compatible with $\{\theta_{ij}\}$ in the sense that
\begin{equazione+}\label{xiv.14}
\theta_j|_i \bydef \theta_j|_{V_{ij}} = {\theta_{ij} \cdot \theta_i|_{V_{ij}}} = {\theta_{ij} \cdot \theta_i|_j}\text. & \text{[mod $\can$]}
\end{equazione+}

For simplicity, let us put $F = {\varphi_!E}$. Our next step will be to define a morphism $\sigma \equiv {\varphi_!\varrho}: {{\s_{\H}}^*F} \to {{\t_{\H}}^*F}$ which is to provide the \H-action on $F$. For each pair $V_i,V_{i'}$ we introduce the abbreviation $\H_{i,i'} = \H(V_i,V_{i'})$. We also write $\H_{ij,i'j'} = \H(V_{ij},V_{i'j'})$. Then the subsets $\H_{i,i'} \subset \mca[1]\H$ form an open cover of the manifold \mca[1]\H. Now, let $g_{i,i'}: \H_{i,i'} \to \G$ be the smooth map obtained by solving the following universal problem
\begin{equazione}\label{xiv.15}
\begin{split}
\xymatrix@C=30pt{\H_{i,i'}\ar[d]_{(\s,\t)}\ar@{-->}[dr]^(.55){g_{i,i'}}\ar@/^1.5pc/[drr]^(.55){\text{~~~~inclusion}} \\ {V_i\times V_{i'}}\ar@/_1.3pc/[dr]_(.3){s_i\times s_{i'}} & \G\ar[d]^{(\s,\t)}\ar[r]^\varphi & \H\ar[d]^{(\s,\t)} \\ & {M\times M}\ar[r]^{f\times f} & {N\times N}\text.\!\!}
\end{split}
\end{equazione}
We can use this map to define a morphism $\sigma_{i,i'}: ({{\s_{\H}}^*F})|_{i,i'} \to ({{\t_{\H}}^*F})|_{i,i'}$ in the category \stack[\H_{i,i'}]{T}
\begin{equazione+}\label{xiv.16'}
\sigma_{i,i'} = {{(\t_{\H}|_{i,i'})^*{\theta_i}^{-1}} \cdot {{g_{i,i'}}^*\varrho} \cdot {(\s_{\H}|_{i,i'})^*\theta_i}}\text. & \text{[mod $\can$]}
\end{equazione+}
By taking into account the equality of mappings
\begin{equazione}\label{xiv.17}
g_{i,i'}|_{j,j'} = {\left({g_{j'i'} \circ \t_{\H}|_{ij,i'j'}}\right) g_{j,j'}|_{i,i'} \left({g_{ji} \circ \s_{\H}|_{ij,i'j'}}\right)}
\end{equazione}
and the identities \refequ{xiv.13'}, \refequ{xiv.14} and \refequ{xiv.16'}, one sees that $\sigma_{i,i'}|_{j,j'} = \sigma_{j,j'}|_{i,i'}$ in \stack[\H_{ij,i'j'}]{T}. Hence the morphisms $\sigma_{i,i'}$ glue together into a unique $\sigma$.

For any morphism $a: R \to R'$ in the category \R{\G}, we obtain a morphism ${\varphi_!a}: {\varphi_!R} \to {\varphi_!R'}$ by setting $b_i = {{s_i}^*a}$ and by observing that
\begin{equazione}\label{xiv.18}
{\theta_{ij}' \cdot b_i|_j} = {b_j|_i \cdot \theta_{ij}} \quad \text{in~}\stack[V_{ij}]T\text.
\end{equazione}
In this way we get a functor of \R{\G} into \R{\H}. The construction of the isomorphisms ${\varphi^* \circ \varphi_!} \simeq \Id_{\R\G}$ and ${\varphi_! \circ \varphi^*} \simeq \Id_{\R\H}$ is left as an exercise.
\end{paragrafo}

\sezione{Smooth Euclidean Fields}\label{v}
\noindent In order to get our reconstruction theory to work effectively, we need to make further hypotheses on the type. We shall say that a smooth tensor stack \stack F is \textit{Euclidean} or, for brevity, that \stack F is a \textit{Euclidean stack,} if it satisfies the following axiomatic conditions (\ref{N.iiiA1}--\ref{N.iiiA6}):

\begin{paragrafo}[Axiom: tensor product and pullback.]\label{N.iiiA1}
\textsl{The canonical natural morphisms \refequ{N.i18} and \refequ{N.i21}
$$%
\left\{\begin{aligned}
& {\sections{E} \otimes_{\smooth[X]} \sections{E'}} \to \sections{({E \otimes E'})}
\\
& f^*(\sections[Y]{F}) \to \sections[X]{({f^*F})}
\end{aligned}\right.
$$%
are surjective (= epimorphisms of sheaves).}

Thus, every local smooth section of ${E \otimes E'}$ will possess, in the vicinity of each point, an expression as a finite linear combination with smooth coefficients of sections of the form ${\zeta \otimes \zeta'}$. Similarly, given any partial smooth section of ${f^* F}$, it will be possible to express it locally as a finite linear combination with coefficients in \smooth[X] of sections of the form ${\eta \circ f}$.

Suppose $E$ is an object of \stack[X]{F}. Let us consider the evaluation map $\sections E(U) \to E_x, \zeta \mapsto \zeta(x)$ defined in \S\ref{iii.3} for a generic open neighbourhood $U$ of the point $x$. When $U$ varies, these maps are evidently mutually compatible, hence on passing to the inductive limit they determine a linear map
\begin{equazione}\label{N.iii1}
(\sections E)_x \to E_x\text, \quad \zeta \mapsto \zeta(x)
\end{equazione}
of the stalk of \sections E at $x$ into the fibre of $E$ at the same point. We call this map the \textit{evaluation} (\textit{of germs}) \textit{at $x$}. It follows from the axiom that for any stack of smooth fields the evaluation of germs at a point is a surjective map. Hence the values $\zeta(x)$ span the fibre $E_x$.
\end{paragrafo}

\begin{paragrafo}[Axiom: criterion for vanishing.]\label{N.iiiA2}
\textsl{Let $a: E \to E'$ be a morphism in \stack[X]F. Suppose that $a_x: E_x \to {E'}_x$ is zero $\forall x \in X$. Then $a = 0$.}

As a first, immediate consequence, one gets that an arbitrary section $\zeta \in \sections E(U)$ vanishes if and only if all the values $\zeta(u)$ are zero as $u$ ranges over $U$. Thus, smooth sections are characterized by their values. Furthermore, by combining this axiom with the former, it follows that the functor $\sections[X]{}: \stack[X]{F} \to \SheavesOfModules{\smooth[\mathnormal X]}$ is faithful.

Each morphism $a: E \to F$ in \stack[X]{F} determines a family of linear maps $\{a_x: E_x \to F_x\}$ and a morphism of sheaves of \smooth[X]-modules $\alpha \equiv \sections a: \sections E \to \sections F$. The link between these two pieces of data is provided by the evaluation maps \refequ{N.iii1}. Namely, for every $x$, the stalk homomorphism $\alpha_x$ and the linear map $a_x$ are ``compatible'': the diagram
\begin{equazione}\label{N.iii5}
\begin{split}
\xymatrix{(\sections E)_x\ar@{->>}[d]_{\text{eval.}}\ar[r]^{\alpha_x} & (\sections F)_x\ar@{->>}[d]^{\text{eval.}} \\ E_x\ar[r]^{a_x} & F_x}
\end{split}
\end{equazione}
commutes. In general, we say that a morphism of sheaves of modules $\alpha: \sections E \to \sections F$ and a family of linear maps $\{a_x: E_x \to F_x\}$ are compatible if (\ref{N.iii5}) commutes for all $x$. Let us call a morphism $\alpha$ of sheaves of modules \textit{representable} if there exists a family of linear maps compatible with $\alpha$.
\end{paragrafo}

\begin{paragrafo}[Axiom: representable morphisms.]\label{N.iiiA3}
\textsl{For each representable morphism $\alpha: \sections E \to \sections F$ there exists an arrow $a: E \to F$ in \stack[X]{F} with $\sections a = \alpha$.}

This axiom will play a role in \refsez{vi}, where we need it in order to construct morphisms of representations by means of fibrewise integration.
\end{paragrafo}

\noindent We say that a form $\phi: {E \otimes E^*} \to \TU$ in the category \stack[X]{F} is a \textit{metric} (\textit{on $E$}) when for every point $x$ the induced form on the fibre $E_x$
\begin{equazione}\label{N.iii6}
{E_x \otimes_\nC {E_x}^*} \to ({E \otimes E^*})_x \xto{\phi_x} \TU_x \can \nC
\end{equazione}
is positive definite Hermitian.

\begin{paragrafo}[Axiom: local metrics.]\label{N.iiiA4}
\textsl{Any object $E$ of the category \stack[X]{F} supports enough local metrics; that is to say, the open subsets $U$ such that one can find a metric on the restriction $E|_U$ cover $X$.}

In general, one can assume only local metrics to exist. Global metrics can be constructed from local ones provided smooth partitions of unity over the manifold $X$ are available.

Let $\phi$ be a metric on $E$. By a \textit{$\phi$-orthonormal frame} (\textit{for $E$}) \textit{about a point $x \in X$} we mean a list of sections $\zeta_1, \ldots, \zeta_d \in \sections E(U)$ defined over a neighbourhood of $x$ such that for all $u \in U$ the vectors $\zeta_1(u), \ldots, \zeta_d(u)$ are orthonormal in $E_u$ and
\begin{equazione}\label{N.iii7}
\mathrm{Span}\{\zeta_1(x),\ldots,\zeta_d(x)\} = E_x\text.
\end{equazione}
We note that orthonormal frames for $E$ exist about each point $x$ at which the fibre $E_x$ is finite dimensional. Indeed, by Axiom \ref{N.iiiA1}, over some neighbourhood $V$ of $x$ one can find local smooth sections $\zeta_1, \ldots, \zeta_d$ with the property that the vectors $\zeta_1(x), \ldots, \zeta_d(x)$ form a basis for $E_x$. Since for all $v \in V$ the vectors $\zeta_1(v), \ldots, \zeta_d(v)$ are linearly dependent if and only if there is a $d$-tuple of complex numbers $(z_1,\ldots,z_d)$ with ${\modulo{z_1}^2 + \cdots + \modulo{z_d}^2} = 1$ and $\txtsum{i=1}{d}{z_i \zeta_i(v)} = 0$, the continuous function
$$%
{V \times \mathbb S^{2d-1}} \to \nR\text, \qquad (v;s_1,t_1,\ldots,s_d,t_d) \mapsto \modulo{\txtsum{k=1}{d}{({s_k + {i t_k}}) \zeta_k(v)}}
$$%
must have a positive minimum at $v=x$, hence a positive lower bound on a suitable neighbourhood $U$ of $x$, so that $\zeta_1(u), \ldots, \zeta_d(u)$ must be linearly independent for all $u \in U$. At this point it is enough to apply the \mbox{Gram}--\mbox{Schmidt} process in order to obtain an orthonormal frame over $U$.

Consider an embedding $e: E' \into E$ of objects of \stack[X]{F}; that is to say, a morphism such that the linear map $e_x: {E'}_x \into E_x$ is injective for all $x$.\footnote{It follows immediately from Axiom \ref{N.iiiA2} that an embedding is a monomorphism. The converse need not be true because the functor $E \mapsto E_x$ does not enjoy any exactness properties. For example, let $a$ be a smooth function on \nR\ such that $a(t) = 0$ if and only if $t = 0$. Then $a$, regarded as an element of $\End(\TU)$, is both mono and epi in $\stack F(\nR)$ while $a_0 = 0: \nC \to \nC$ is neither injective nor surjective.} Suppose there exists a global metric $\phi$ on the object $E$. Also assume that $E'$ is a locally trivial object of locally finite rank. Then $e$ admits a cosection, i.e.\ there exists a morphism $p: E \to E'$ with ${p \circ e} = \id$. To prove this, note first of all that the metric $\phi$ induces a metric $\phi'$ on $E'$. For each point $x$ there exists a $\phi'$-orthonormal frame $\zeta'_1, \ldots, \zeta'_d \in [\sections{E'}](U)$ for $E'$ about $x$, since ${E'}_x$ is finite dimensional. Let $\zeta^i$ be the composite
\begin{equazione}\label{N.iii10}
E|_U \can {E|_U \otimes \TU_U} \can {E|_U \otimes {\TU|_U}^*} \xto{\;E|_U \otimes {\zeta_i}^*\;} {E|_U \otimes {E|_U}^*} \xto{\quad\phi|_U\quad} \TU_U\text,
\end{equazione}
where $\zeta_i \equiv [\sections e(U)](\zeta'_i)$. Define $p_U$ as ${({\zeta'_1 \oplus \cdots \oplus \zeta'_d}) \cdot ({\zeta^1 \oplus \cdots \oplus \zeta^d})}$ (``orthogonal projection onto $E'|_U$''). Our claim follows from Axiom \ref{N.iiiA2}.

By using the last remark, and once more the existence of local orthonormal frames, one can show that if the dimension of the fibres of an object $E$ of \stack[X]{F} is finite and locally constant over $X$ then $E$ is locally trivial of locally finite rank.
\end{paragrafo}

\begin{lemma}\label{N.iii4}
(Let \stack F be a Euclidean stack.) Let $X$ be a paracompact manifold and $i_S: S \into X$ a closed submanifold. Let $E, F$ be objects of \stack[X]{F} and suppose that $E' \equiv E|_S$ is locally free of locally finite rank over $S$. Put $F' = F|_S$. Then every morphism $a': E' \to F'$ in \stack[S]{F} can be extended to a morphism $a: E \to F$ in \stack[X]{F}.
\end{lemma}
\begin{proof}
Fix a point $s \in S$. There exists an open neighbourhood $A$ of $s$ in $S$ such that there is a trivialization $E'|_A \iso {\TU_A \oplus \cdots \oplus \TU_A}$ over $A$. Let $\zeta'_1, \ldots, \zeta'_d \in \sections{E'}(A)$ be the corresponding frame of local sections. Also, let $U$ be any open subset of $X$ such that ${U \cap S} = A$. After taking $U$ and $A$ smaller about $s$ if necessary, it is no loss of generality to assume, by Axiom \refcnt{N.iiiA1}, that there are local sections $\zeta_1, \ldots, \zeta_d \in \sections E(U)$ with $\zeta'_k = {\zeta_k \circ i_S}$, $k = 1,\ldots,d$. The values $\zeta_k(x)$, $k = 1,\ldots,d$ must be linearly independent in the fibre $E_x$ because the same is true of the values $\zeta'_k(s)$, $k = 1,\ldots,d$ in ${E'}_s$. This implies that if $U$ is small enough then the morphism $\zeta \equiv {\zeta_1 \oplus \cdots \oplus \zeta_d}: {\TU_U \oplus \cdots \oplus \TU_U} \to E|_U$ is an embedding and admits a cosection $p: E|_U \to {\TU_U \oplus \cdots \oplus \TU_U}$, as observed above.

Set $\eta'_k = [\sections{a'}(A)](\zeta'_k) \in [\sections{F'}](A)$. As before, it is no loss of generality to assume that there are sections $\eta_1,\ldots,\eta_d \in \sections F(U)$ with $\eta'_k = {\eta_k \circ i_S}$. These can be combined into a morphism $\eta: {\TU_U \oplus \cdots \oplus \TU_U} \to F|_U$ ($d$\nobreakdash-fold direct sum). Then one can take the composition
\begin{equazione}
E|_U \xto{\:p\:} {\TU_U \oplus \cdots \oplus \TU_U} \xto{\:\eta\:} F|_U\text.
\end{equazione}
It is immediate to check that the restriction of this morphism to the submanifold $A \into U$ coincides with $a'|_A$, up to the canonical identifications $(E|_U)|_A \can E'|_A$ and $(F|_U)|_A \can F'|_A$.

One concludes the proof by using a partition of unity over $X$.
\end{proof}

\begin{paragrafo}[Axiom: dimension.]\label{N.iiiA5}
\textsl{It is required of the canonical ``pseudo-''tensor* functor $\stack[\pt]{F} \longto \ComplexVectorSpaces$ (\ref{N.i10}) that
\begin{elenco}
\item it is fully faithful;
\item it factors through the subcategory whose objects are the finite dimensional vector spaces, in other words the vector space $E_*$ (\ref{N.i10}) is finite dimensional for all $E$ in \stack[\pt]{F};
\item it is a genuine tensor* functor, i.e.\ (\ref{N.i7}) and (\ref{N.i8}) are isomorphisms of sheaves for $X = \pt$.
\end{elenco}}

It follows from this axiom that the functor $E \mapsto E_x$ is a true tensor* functor (in general it is only a ``pseudo-''tensor* functor).
\end{paragrafo}

\noindent We shall say that an object $E$ in \stack[X]{F} is \textit{locally finite} if the sheaf \sections{E} is a locally finitely generated \smooth[X]-module; in other words, $E$ is locally finite if the manifold $X$ admits a cover by open subsets $U$ such that for each of them there is an epimorphisms of sheaves of modules
\begin{equazione}\label{N.iii9}
{\smooth[U] \oplus \cdots \oplus \smooth[U]} \xto{\quad\text{epi}\quad} (\sections{E})|_U\text.
\end{equazione}

\begin{paragrafo}[Axiom: local finiteness.]\label{N.iiiA6}
\textsl{For every manifold $X$, all the objects of the category \stack[X]{F} are locally finite.}
\end{paragrafo}

\begin{paragrafo}[Example: smooth Hilbert fields.]\label{xvi.1}
By a ``smooth Hilbert field'' we mean an object \f[*]H consisting of a family $\{H_x\}$ of complex Hilbert spaces indexed by the set of points of a manifold $X$ and a sheaf \f{H} of \smooth[X]-modules of local sections of $\{H_x\}$ subject to the following conditions:
\begin{elenco}
\item $\{\zeta(x): \zeta \in (\f{H})_x\}$, where $(\f{H})_x$ indicates the stalk at $x$, is a dense linear subspace of $H_x$;
\item for each open subset $U$, and for all sections $\zeta, \zeta' \in \f{H}(U)$, the function \scalare{\zeta}{\zeta'} on $U$ defined by $u \mapsto \scalare{\zeta(u)}{\zeta'(u)}$ is smooth.
\end{elenco}
We refer to the manifold $X$ as the ``base'' of \f[*]H; we will also say that \f[*]H is a smooth Hilbert field over $X$. Let \f[*]H and \f[*]{H'} be smooth Hilbert fields over a manifold $X$. A morphism of \f[*]H into \f[*]{H'} is a family of bounded linear maps $\{a_x: H_x \to {H'}_x\}$ indexed by the set of points of $X$ such that for each open subset $U$ and for each $\zeta \in \f{H}(U)$ the section over $U$ given by $u \mapsto {a_u \cdot \zeta(u)}$ belongs to $\f{H'}(U)$. Smooth Hilbert fields over $X$ and their morphisms form a category which we shall denote by \HB{X}.

\sloppy
Suppose \f[*]H and \f[*]G are Hilbert fields over a manifold $X$. Consider the bundle of tensor products $\{{H_x \otimes G_x}\}$. For any pair of sections $\zeta \in [\f{H}](U)$ and $\eta \in [\f{G}](U)$ we let ${\zeta \otimes \eta}$ denote the section of the bundle $\{{H_x \otimes G_x}\}$ defined over $U$ by $u \mapsto {\zeta(u) \otimes \eta(u)}$. The correspondence
\begin{equazione}\label{N.iiiH3}
U\: \mapsto\: {\smooth(U) \bigl\{{\zeta \otimes \eta}: \zeta \in [\f{H}](U), \eta \in [\f{G}](U)\bigr\}}
\end{equazione}
defines a sub-presheaf of the sheaf of sections of $\{{H_x \otimes G_x}\}$. [Here ${\smooth(U) \{\cdots\}}$ stands for the $\smooth(U)$-module spanned by $\{\cdots\}$.] Let ${\f[*]H \otimes \f[*]G}$ denote the Hilbert field over $X$ given by the bundle $\{{H_x \otimes G_x}\}$ together with the sheaf of sections generated by the presheaf \refequ{N.iiiH3}. We call ${\f[*]H \otimes \f[*]G}$ the tensor product of \f[*]H and \f[*]G. Observe that for all morphisms $\f[*]H \xto\alpha \f[*]{H'}$ and $\f[*]G \xto\beta \f[*]{G'}$ of Hilbert fields over $X$ the bundle of bounded linear maps $\{{a_x \otimes b_x}\}$ yields a morphism ${\alpha \otimes \beta}$ of ${\f[*]H \otimes \f[*]G}$ into ${\f[*]{H'} \otimes \f[*]{G'}}$.

\fussy
One gets the conjugate field $\f[*]H^*$ of a Hilbert field \f[*]H by taking the bundle $\{{H_x}^*\}$ of conjugate spaces along with the local sections of \f[*]H regarded as local sections of $\{{H_x}^*\}$.

With the obvious tensor unit and the obvious \textit{ACU} constraints, these operations turn \HB{X} into a tensor* category. It remains to define a tensor* functor $f^*: \HB{Y} \to \HB{X}$ for each smooth map $f: X \to Y$, along with suitable fibred tensor category constraints.

Let \f[*]G be a Hilbert field over $Y$. The pull-back of \f[*]G along $f$, to be denoted by ${f^* \f[*]G}$, is the smooth Hilbert field over $X$ whose associated bundle of Hilbert spaces is $\{G_{f(x)}\}$ and whose associated sheaf of sections is generated by the following presheaf of sections of the bundle $\{G_{f(x)}\}$:
\begin{equazione}\label{N.iiiH4}
U\: \mapsto\: {\smooth[X](U) \bigl\{{\eta \circ f}: \eta \in [\f{G}](V), V \supset f(U)\bigr\}}\text.
\end{equazione}
For every morphism $\beta: \f[*]G \to \f[*]{G'}$ of Hilbert fields over $Y$, the family of bounded linear maps $\{b_{f(x)}\}$ defines a morphism ${f^* \beta}: {f^* \f[*]G} \to {f^* \f[*]{G'}}$ of Hilbert fields over $X$. The operation $\f[*]G \mapsto {f^* \f[*]G}$ defines a ``strict'' tensor* functor of \HB{Y} into \HB{X}, in other words one has the identities
$$%
f^*({\f[*]G \otimes \f[*]{G'}}) = {{f^* \f[*]G} \otimes {f^* \f[*]{G'}}}\text, \quad f^*(\TU_Y) = \TU_X \quad \text{and} \quad {f^* \f[*]{(G^*)}} = ({f^* \f[*]G})^*\text.
$$%
Finally, the identities of tensor* functors
$$%
({g \circ f})^* = {f^* \circ g^*} \quad \text{and} \quad {\id_X}^* = \Id
$$%
provide the required fibred tensor category constraints.

The fibred tensor category $X \mapsto \HB{X}$ is a smooth tensor stack satisfying Axioms \refcnt{N.iiiA1}, \refcnt{N.iiiA2} and \refcnt{N.iiiA4}. However, as it does not satisfy the other axioms, it is not an example of a Euclidean stack.
\end{paragrafo}

\begin{paragrafo}[Example: smooth Euclidean fields.]\label{xvi.2}
Let \EB{X} denote the full subcategory of \HB{X} consisting of all \f[*]E whose associated sheaf of sections \f{E} is locally finite. We refer to the objects of this subcategory as ``smooth Euclidean fields'' (over $X$).

Observe that \EB{X} is a tensor* subcategory of \HB{X}. Indeed, since the smooth tensor stack of smooth Hilbert fields satisfies Axiom \refcnt{N.iiiA1}, the locally finite \smooth[X]-module ${\f{E} \otimes_{\smooth[X]} \f{E'}}$ surjects onto the \smooth[X]-module $\f{({E \otimes E'})}$. For similar reasons, for any map $f: X \to Y$ the pullback functor $f^*: \HB{Y} \to \HB{X}$ must carry \EB{Y} into \EB{X}.

The smooth tensor stack $X \mapsto \EB{X}$ also satisfies Axioms \refcnt{N.iiiA3}, \refcnt{N.iiiA5} and \refcnt{N.iiiA6} and is therefore Euclidean.
\end{paragrafo}

\sezione{Construction of Equivariant Maps}\label{vi}
Let \stack F denote an arbitrary Euclidean stack. \stack F is to be regarded as fixed throughout the entire section.

\begin{lemma}\label{N.iv2}
Let \G\ be a (locally) transitive Lie groupoid, and let $X$ be its base manifold. Take an arbitrary representation $(E,\varrho) \in \R[F]{\G}$. Then $E$ is locally trivial in \stack[X]{F}.
\end{lemma}
\begin{proof}
Local transitivity means that the mapping $(\s,\t): \G \to {X\times X}$ is a submersion. Fix a point $x \in X$. Since $(x,x)$ lies in the image of the map $(\s,\t)$, the latter admits a local smooth section ${U\times U} \to \G$ over some open neighbourhood of $(x,x)$. Let us consider the restriction $g: U \to \G$ of this section to $U = {U\times \{x\}}$.

Let $x: \pt \to X$ denote the map $\pt \mapsto x$. By Axiom \refcnt{N.iiiA5}, there is a trivialization for ${x^*E}$ in \stack[\pt]{F}. We pull $\varrho$ back to $U$ along the smooth map $g$, and observe that there is a unique factorization of ${\t\circ g}$ through \pt\ (collapse $c: U \to \pt$ followed by $x: \pt \to X$). Since $\varrho$ is an isomorphism,
\begin{multline}\notag%
E|_U = {{i_U}^*E} = {({\s\circ g})^*E} \can {g^*\s^*E} \xto{\:\:g^*\varrho\:\:} {g^*\t^*E} \can {({\t\circ g})^*E} =
\\[\medskipamount]
= {({x\circ c})^*E} \can {c^*(x^*E)} \iso {c^*({\TU \oplus \cdots \oplus \TU})} = {\TU_U \oplus \cdots \oplus \TU_U}
\end{multline}%
provides a trivialization for $E|_U$ in \stack[U]{F}.
\end{proof}
Let $i_S: S \into X$ be an invariant immersed submanifold. The pullback $\G|_S$ of $\G$ along $i_S$ is well-defined and is a Lie subgroupoid of \G.\footnote{In general, a Lie subgroupoid is a Lie groupoid homomorphism $(\varphi,f)$ such that both $\varphi$ and $f$ are injective immersions; compare for instance \cite{MoeMrc06}.} [Observe that $\G|_S = \G^S = {\s_\G}^{-1}(S)$.] In the special case of an orbit immersion, $\G|_S$ will be transitive over $S$. Then the lemma says that for any $(E,\varrho) \in {\Ob\,\R\G}$ the restriction $E|_S$ is a locally trivial object of \stack[S]{F}.

\begin{paragrafo}[Alternative description of representations.]\label{N.iv1}%
The notion of representation with which we have been working so far is completely intrinsic. We were able to prove all results by means of purely formal arguments involving only manipulations of commutative diagrams. For the purposes of the present section, however, we have to change our point of view.

Let \G\ be a Lie groupoid. Consider a representation $\s^*E \xto\varrho \t^*E$ of \G. Each arrow $g$ determines a linear map $\varrho(g): E_{\s(g)} \to E_{\t(g)}$ defined via the commutativity of the diagram
\begin{equazione}\label{N.iv1a}
\begin{split}
\xymatrix{[g^*\s^*E]_*\ar[d]^{[g^*\varrho]_*}\ar[r]^-{[\can]_*} & [\s(g)^*E]_*\ar@{=}[r]^-{\text{def}} & E_{\s(g)}\ar@{-->}[d]^{\varrho(g)} \\ [g^*\t^*E]_*\ar[r]^-{[\can]_*} & [\t(g)^*E]_*\ar@{=}[r]^-{\text{def}} & E_{\t(g)}}
\end{split}
\end{equazione}
[the notation \refequ{N.i10} is in use]. It is routine to check that the conditions \textrm{i)} and \textrm{ii)} in the definition of a representation (beginning of \refsez{iv}) imply that the correspondence $g \mapsto \varrho(g)$ is multiplicative, i.e.\ that $\varrho(g'g) = {\varrho(g') \circ \varrho(g)}$ and $\varrho(x) = \id$ for all points $x$ of the base manifold $X$.

Fix an arbitrary arrow $g_0$. Let $\zeta \in \sections E(U)$ be a section, defined over a neighbourhood of $\s(g_0)$ in $X$. Recall that, according to \refequ{N.i19}, $\zeta$ will determine a section ${\zeta\circ\s} \in \sections[{\G}]{(\s^*E)}(\G^U)$ at which the morphism of sheaves of modules \sections{\varrho} can be evaluated so as to get a section of $\t^*E$ over $\G^U$. Now, Axiom \refcnt{N.iiiA1} implies that there exists an open neighbourhood $\Gamma$ of $g_0$ in $\G^U$ over which the latter section can be expressed as a finite linear combination with coefficients in $\C^\infty(\Gamma)$ of sections of the form ${\zeta'_i\circ\t}$, with $\zeta'_i$ ($i = 1, \cdots, d$) defined over $\t(\Gamma)$. In symbols,
\begin{equazione}\label{N.iv1b}%
\bigl[\sections\varrho\,(\Gamma)\bigr]({\zeta \circ \s|_\Gamma}) = \txtsum{i=1}d{r_i ({\zeta'_i\circ\t})|_\Gamma}\text,
\end{equazione}%
with $r_1, \ldots, r_d \in \C^\infty(\Gamma)$ and $\zeta'_1, \ldots, \zeta'_d \in [\sections E][\t(\Gamma)]$. This equality can be ``evaluated'' at $g \in \Gamma$ to get
\begin{equazione}\label{N.iv1c}%
{\varrho(g) \cdot \zeta(\s[g])} = \txtsum{i=1}d{r_i(g) \zeta'_i(\t[g])}\text.
\end{equazione}%
By Axiom \refcnt{N.iiiA3}, any multiplicative operation $g \mapsto \varrho(g)$, locally of the form \refequ{N.iv1c}, comes from a representation of \G\ on $E$.
\end{paragrafo}

\begin{paragrafo}[Preliminary extension.]\label{vi.1}
Suppose \G\ proper hereafter. Fix a point $x_0 \in X$, and let $G_0$ denote the isotropy group at $x_0$. It is evident from \refequ{N.iv1c} that
\begin{equazione}\label{N.iv3}%
\varrho_0: G_0 \to \GL(E_0)\text, \quad g \mapsto \varrho(g)
\end{equazione}%
is a continuous representation of the compact Lie group $G_0$ on the finite dimensional vector space $E_0$ (the fibre of $E$ at $x_0$).

Suppose another \G-action $(F,\sigma)$ is given, along with some $G_0$-equivariant linear map $A_0: E_0 \to F_0$. Let $S_0 \into X$ be the orbit through $x_0$. Our remarks about Morita equivalences in \S\ref{iv} say there exists a unique morphism $A': (E',\varrho') \to (F',\sigma')$ in \R{\G'} [the primes here signify that we are taking the corresponding restrictions to $S_0$] such that $(A')_0 = A_0$. In fact, for every point $z \in S_0$ and arrow $g \in \G(x_0,z)$, one has
\begin{equazione}\label{N.iv4}%
(A')_z = {\sigma(g) \cdot A_0 \cdot {\varrho(g)}^{-1}}: E_z \to F_z\text.
\end{equazione}%
By Lemma \ref{N.iv2}, $E'$ is a locally trivial object of \stack[S_0]{F}. Then Lemma \ref{N.iii4} yields a global morphism $a: E \to F$ extending $A'$ and hence, a fortiori, $A_0$. We proceed to ``average out'' this $a$ to make it \G-equivariant, as follows.
\end{paragrafo}

\begin{paragrafo}[Averaging operators.]\label{vi.2}
Fix an arbitrary (right invariant, normalized) Haar system $\mu = \{\mu^x\}$ on the (proper) Lie groupoid \G. We shall construct, for each pair of \G-actions $R = (E,\varrho)$ and $S = (F,\sigma)$, a linear operator
\begin{equazione}\label{N.ivB1}%
\Av{} = \Av\mu: \Hom_{\stack[M]F}(E,F) \to \Hom_{\R\G}(R,S)
\end{equazione}%
(averaging operator), with the property that $\Av{}(a) = a$ whenever $a$ already belongs to $\Hom_{\R\G}(R,S)$. More generally, if $S$ is an invariant submanifold over which $a$ restricts to an equivariant morphism, $\Av{}(a)|_S = a|_S$.

We start from a very simple remark. Suppose sections $\zeta \in \sections E(U)$ and $\eta_1, \ldots, \eta_n \in \sections F(U)$ are given such that $\eta_1, \ldots, \eta_n$ are local generators for \sections F over $U$. Then for each $g_0 \in \G^U$ there exists an open neighbourhood $\Gamma \subset \G^U$ of $g_0$ along with smooth functions $\phi_1, \ldots, \phi_n$ on $\Gamma$ such that
\begin{equazione}\label{N.ivB3}%
{\sigma(g)^{-1} \cdot a_{\t(g)} \cdot \varrho(g) \cdot \zeta(\s[g])} = \txtsum{j=1}n{\phi_j(g)\eta_j(\s[g])}
\end{equazione}%
for all $g \in \Gamma$. To see this, note that\inciso{as observed in \refequ{N.iv1c}}there are an open neighbourhood $\Gamma$ of $g_0$ in $\G^U$ and local smooth sections $\zeta'_1, \ldots, \zeta'_m$ of $E$ over $U' = \t(\Gamma)$ such that ${\varrho(g) \zeta(\s[g])} = \txtsum{i=1}m{r_i(g) \zeta'_i(\t[g])}$ for some functions $r_1, \ldots, r_m \in \C^\infty(\Gamma)$. For $i = 1, \ldots, m$, put $\eta'_i = {\sections a(U')(\zeta'_i)} \in \sections F(U')$. Since $\Gamma^{-1}$ is a neighbourhood of ${g_0}^{-1}$, by using the hypothesis that the $\eta_j$'s are generators we can also assume $\Gamma$ to be so small that for each $i = 1, \ldots, m$ there exist $s_{1,i}, \ldots, s_{n,i} \in \C^\infty(\Gamma^{-1})$ with ${\sigma(g^{-1}) \eta'_i(\t[g])} = \txtsum{j=1}n{s_{j,i}(g^{-1}) \eta_j(\s[g])}$ for each $g \in \Gamma$. This proves \refequ{N.ivB3}.

Put $\alpha = \sections a$. We can use the last remark to obtain a new morphism of sheaves of modules $\tilde\alpha: \sections E \to \sections F$, as follows. Let $\zeta$ be a local section of $E$ defined over an open subset $U$ so small that by Axiom \ref{N.iiiA6} there exists a system $\eta_1, \ldots, \eta_n$ of local generators for \sections F over $U$. For each $g_0 \in \G^U$, select an open neighbourhood $\Gamma(g_0)$ along with smooth functions $\phi^{g_0}_1, \ldots, \phi^{g_0}_n \in \C^\infty\bigl(\Gamma(g_0)\bigr)$ satisfying \refequ{N.ivB3}. Then choose a smooth partition of unity over $\G^U$ $\{\theta_i: i \in I\}$ subordinated to the $\Gamma(g_0)$, and put
\begin{equazione}\label{N.ivB4}%
{\tilde\alpha(U)}(\zeta) = \txtsum{j=1}n{\Phi_j \eta_j} \quad \text{where} \quad \Phi_j(u) = {\int_{\G^u} \txtsum{i\in I}{}{\theta_i(g) {\phi^i}_j(g)}\:\mathit d\mu^u(g)}\text.
\end{equazione}%
Some arbitrary choices are involved here, so one has to make sure that this is a good definition. If we look at \refequ{N.ivB3} for $x = \s(g)$ fixed, we recognize that the operation $g \mapsto {\sigma(g)^{-1} \cdot a_{\t(g)} \cdot \varrho(g) \cdot \zeta(x)}$ defines a smooth mapping on the manifold $\G^x$ with values in the finite dimensional vector space $F_x$. For each $v \in E_x$, there is some local section $\zeta$ about $x$ such that $\zeta(x) = v$, so one is allowed to take the integral
\begin{equazione}\label{N.ivB6}%
\varkappa_x(v) = {\int_{\G^x} {\sigma(g)^{-1} \cdot a_{\t(g)} \cdot \varrho(g) \cdot v}\:\mathit d\mu^x(g)}\text.
\end{equazione}%
This defines, for each base point $x$, a linear map $\varkappa_x: E_x \to F_x$. Now,
\begin{align*}
{[\tilde\alpha(U)(\zeta)]}(u) &= \txtsum{j=1}n{\Phi_j(u) \eta_j(u)} = \txtsum{j=1}n{{\int_{\G^u} \txtsum{i\in I}{}{\theta_i(g) {\phi^i}_j(g)}\,\mathit d\mu^u(g)}\: \eta_j(u)}
\\
&= {\int_{\G^u} \txtsum{i\in I}{}{\theta_i(g) \txtsum{j=1}n{{\phi^i}_j(g) \eta_j(\s[g])}}\:\mathit d\mu^u(g)}
\\
&= {\int_{\G^u} \txtsum{i\in I}{}{\theta_i(g)} {\Bigl[\sigma(g)^{-1} \cdot a_{\t(g)} \cdot \varrho(g) \cdot \zeta(\s[g])\Bigr]}\,\mathit d\mu^u(g)}
\\
&= {\varkappa_u \cdot \zeta(u)}\text.
\end{align*}
It follows from Axiom \ref{N.iiiA2} that the section $\tilde\alpha(U)(\zeta)$ in \refequ{N.ivB4} does not depend on any of the auxiliary choices we made in order to define it (as the $\varkappa_u$ don't).

We define $\Av{}(a)$ as the unique morphism $\tilde a: E \to F$ with $\sections{(\tilde a)} = \tilde\alpha$. [Its existence follows from the preceding computation and Axiom \ref{N.iiiA3}, its uniqueness from Axiom \ref{N.iiiA2}.] It remains to show that $\Av{\mu}$ is a projection operator onto $\Hom_{\R\G}(R,S)$. We will leave the verification to the reader.
\end{paragrafo}

\noindent Summing up \ref{vi.1} and \ref{vi.2}, one gets

\begin{proposizione}\label{N.iv5}
Suppose \G\ is proper, and let $x_0$ be a base point. For each pair of \G-actions $R = (E,\varrho)$ and $S = (F,\sigma)$, and for each $G_0$-equivariant linear map $A_0: E_0 \to F_0$, there exists in \R{\G} a morphism $a: R \to S$ extending $A_0$.
\end{proposizione}

\noindent By applying the averaging operator to a randomly chosen Hermitian metric, we get the existence of invariant metrics

\begin{proposizione}\label{N.iv6}
Let $R = (E,\varrho)$ be a representation of a proper Lie groupoid \G. Then there exists a \G-invariant metric on $E$, that is, a metric on $E$ which is at the same time a morphism ${R\otimes R^*} \to \TU$ in \R\G.
\end{proposizione}

\noindent By a $\varrho$-invariant partial section of $E$ over an invariant submanifold $S$ of the base of \G\ we mean a section of $E|_S$ over $S$ which is at the same time a morphism in $\R{\G|_{\mathnormal S}}$. Lemma \ref{N.iii4} in combination with \ref{vi.2} yields

\begin{proposizione}\label{iv7}
Let $S$ be a closed invariant submanifold of the base of a proper Lie groupoid \G. Let $R = (E,\varrho)$ be a representation of \G. Then each $\varrho$-invariant partial section of $E$ over $S$ can be extended to a global $\varrho$-invariant section of $E$.
\end{proposizione}

\noindent A function $\varphi$ defined on an arbitrary subset $S$ of a manifold $X$ is called smooth when for each $x \in X$ one can find an open neighbourhood $U$ of $x$ in $X$ and a smooth function on $U$ which agrees with $\varphi$ on ${U\cap S}$.

\begin{proposizione}\label{iv8}
Let $S$ be an invariant subset of the base manifold $X$ of a proper Lie groupoid \G. Suppose $\varphi$ is a smooth invariant (viz.\ constant along the \G-orbits) function on $S$. Then there exists a smooth invariant function extending $\varphi$ on all of $X$.
\end{proposizione}
\begin{proof}
Average out any smooth extension of $\varphi$ obtained by means of a partition of unity over $X$.
\end{proof}

\sezione{$\C^\infty$ Fibre Functors}\label{vii}
We keep on working with a generic Euclidean stack \stack F. Let $M$ be a paracompact smooth manifold.

\begin{definizione}\label{xviii.1}
By a \textit{fibre functor} (\textit{of type \stack F}) \textit{over $M$}, or \textit{with base $M$}, we mean a faithful tensor* functor
\begin{equazione}\label{N.v0}
\fifu: \Kt \longto \stack[M]{F}
\end{equazione}
defined on some tensor* category \Kt.
\end{definizione}
When a fibre functor \fifu\ is assigned over $M$, one can construct a groupoid \tannakian{\fifu} having the points of $M$ as base points. Under reasonable assumptions, it is possible to endow \tannakian{\fifu} with a natural structure of topological groupoid; the choice of a topology is dictated by the idea that the objects of \Kt\ should give rise to continuous representations of \tannakian{\fifu} and that, vice versa, continuity of these representations should be enough to characterize the topology. An improvement of the same idea leads one to introduce a certain $\C^\infty$ functional structure on the space of arrows of \tannakian{\fifu}. (Recall \S\ref{C^infty-Sp} and \S\ref{C^infty-Gpd}.) When \tannakian{\fifu} is a $\C^\infty$-groupoid relative to this particular $\C^\infty$-structure, we say that \fifu\ is a \textit{$\C^\infty$ fibre functor.} In detail, these constructions read as follows.

\begin{paragrafo}[The Tannakian groupoid \tannakian{\fifu}.]\label{vii.1}
Let $x$ be a point of $M$; the same symbol will be used to denote the corresponding (smooth) map $\pt \to M$. Consider the tensor* functor (cfr \S\refcnt{iii.3} and \S\refcnt{N.iiiA5})
\begin{equazione}\label{N.v1}
\stack[M]{F} \longto \VectorSpaces\text, \quad E \mapsto E_x\text.
\end{equazione}
Let $\fifu_x$ denote the composite tensor* functor
\begin{equazione}\label{N.v2}
\Kt \xto{\:\:\fifu\:\:} \stack[M]{F} \xto{\:(\text-)_x\:} \VectorSpaces\text, \quad R \mapsto (\fifu(R))_x\text.
\end{equazione}

We define two groupoids \tannakian{\fifu;\nC} and \tannakian{\fifu;\nR} over $M$ by putting
\begin{equazione}\label{N.v3}
\tannakian{\fifu;\nC}(x,x') = \Iso^\otimes(\fifu_x,\fifu_{x'}) \quad \text{and} \quad \tannakian{\fifu;\nR}(x,x') = \Iso^{\otimes,*}(\fifu_x,\fifu_{x'})
\end{equazione}
where $x, x' \in M$. (Recall that the right-hand term in the second equality denotes the set of all self-conjugate tensor preserving natural isomorphisms.) By setting $({\lambda' \cdot \lambda})(R) = {\lambda'(R) \circ \lambda(R)}$ and $x(R) = \id$, in each case we obtain a structure of groupoid over $M$. The relationship between \tannakian{\fifu;\nC} and its subgroupoid \tannakian{\fifu;\nR} can be clarified by introducing the ``conjugation involution'' of \tannakian{\fifu;\nC}: this sends an arrow $\lambda$ to the arrow $\overline\lambda$ defined by setting $\overline\lambda(R) = \lambda(R^*)^*$ [up to $\can$]. The elements of \tannakian{\fifu;\nR} are the fixed points of the conjugation involution.

The groupoid \tannakian{\fifu;\nR} shall be referred to as the \textit{Tannakian groupoid} (\textit{associated with \fifu}). We will abbreviate \tannakian{\fifu;\nR} into \tannakian{\fifu}.
\end{paragrafo}

\begin{paragrafo}[Representative functions.]\label{vii.2}
Let $R \in \Ob(\Kt)$ be arbitrary and let $\phi$ be any metric on $\fifu(R)$. For each pair of global sections $\zeta, \zeta' \in \sections{({\fifu R})}(M)$ we introduce the function
\begin{multiriga}\label{N.v8}
r_{R,\phi,\zeta,\zeta'}: \tannakian{\fifu} \to \nC\text,\quad \lambda \mapsto \smash{\bigsca{\lambda(R) \cdot \zeta(\s[\lambda])}{\zeta'(\t[\lambda])}}_\phi
\\
\bydef \phi_{\t(\lambda)}\bigl({\lambda(R) \cdot \zeta(\s[\lambda])},\zeta'(\t[\lambda])\bigr)\text.
\end{multiriga}
We put
\begin{equazione}\label{N.v9}
\Rf = \{r_{R,\phi,\zeta,\zeta'}: R \in \Ob(\Kt)\text, \text{~$\phi$~metric~on~$\fifu(R)$,~} \zeta, \zeta' \in \sections{({\fifu R})}(M)\}\text.
\end{equazione}
We call the elements of \Rf\ \textit{representative functions.} Observe that \Rf\ is a complex algebra of functions on \tannakian{\fifu}, closed under the operation of taking the complex conjugate. This implies that the real and imaginary parts of any function of \Rf\ also belong to \Rf. Thus, if we let $\nR[\Rf] \subset \Rf$ denote the subset of all real valued functions, we have $\Rf = {\nC\otimes{\nR[\Rf]}}$.
\end{paragrafo}

\begin{paragrafo}[Topology and $\C^\infty$-structure.]\label{vii.3}
We endow \tannakian{\fifu} with the smallest topology making all representative functions continuous. As a consequence of the existence of metrics on any object of \stack[M]F, the topological space \tannakian{\fifu} is necessarily Hausdorff. The functions in $\nR[\Rf]$ generate a functional structure on the space \tannakian{\fifu}. One can complete this functional structure to a $\C^\infty$-structure $\Rf^\infty$ as explained in \S\ref{C^infty-Sp}.

We remark that the source map of the groupoid \tannakian{\fifu} is a morphism of $\C^\infty$-spaces relative to the $\C^\infty$-structure $\Rf^\infty$. The same statement is true of the target map and the unit section. However, without stronger assumptions on the fibre functor \fifu\ we are at present unable to show that \tannakian{\fifu} is a $\C^\infty$-groupoid relative to $\Rf^\infty$. It might be the case that not every fibre functor is $\C^\infty$. We will see later on that the standard forgetful functor associated with a ``reflexive'' groupoid is always a $\C^\infty$ fibre functor. This is in fact the only case of interest in connection with the proof of our reconstruction theorem.
\end{paragrafo}

\begin{paragrafo}[Invariant metrics.]\label{vii.4}
Let $R \in \Ob(\Kt)$. We say that a metric $\phi$ on $\fifu(R)$ is \textit{\fifu\nobreakdash-invariant} if there is a Hermitian form $m: {R\otimes R^*} \to \TU$ such that $\phi$ coincides with the induced form
\begin{equazione}\label{N.v12}
{\fifu(R) \otimes \fifu(R)^*} \can \fifu({R\otimes R^*}) \xto{\fifu(m)} \fifu(\TU) \can \TU\text.
\end{equazione}
Note that, by the faithfulness of \fifu, there is at most one such $m$.
\end{paragrafo}

\begin{definizione}\label{xix.1}
A fibre functor $\fifu: \Kt \longto \stack[M]{F}$ will be called \textit{proper} if
\begin{elenco}
\item the continuous mapping $(\s,\t): \tannakian{\fifu} \to {M\times M}$ is proper, and
\item for every $R \in \Ob(\Kt)$, the object $\fifu(R)$ supports an \fifu-invariant metric.
\end{elenco}
We can express the second condition more succinctly by saying that ``{there are enough \fifu-invariant metrics}''.
\end{definizione}

\begin{paragrafo}[Example.]\label{N.v13}
As an example of a proper fibre functor, we mention [recall \S\ref{iv}] the standard forgetful functor $\forgetfulfunctor{\G}: \R{\G} \to \stack[M]{F}$ associated with the representations of type \stack F of a proper Lie groupoid \G\ over $M$.

To begin with, we observe that there is a homomorphism of groupoids
\begin{equazione}\label{N.v14}
\G \longto \tannakian{\forgetfulfunctor{\G}}
\end{equazione}
which sends $g$ to the natural transformation assigning each object $(E,\varrho)$ of the category \R{\G} the isomorphism $\varrho(g)$ [cfr \S\ref{N.iv1}]. This homomorphism is evidently a morphism of $\C^\infty$-spaces and, in particular, a continuous map. It will be established in the next section that (\ref{N.v14}) is a surjection. The properness of $(\s,\t): \tannakian{\forgetfulfunctor{\G}} \to {M\times M}$ is then an immediate consequence of the properness of $(\s,\t): \G \to {M\times M}$. The existence of enough invariant metrics has been proved in the preceding section (Proposition \ref{N.iv6}).
\end{paragrafo}

Let $\Rf' \subset \Rf$ be the set of all representative functions of the form $r_{R,\phi,\zeta,\zeta'}$ where $\phi$ is an \fifu-invariant metric. Note that $\Rf'$ is a subalgebra of $\Rf$ closed under complex conjugation.
\begin{lemma}\label{N.v15}
Let \fifu\ be a proper fibre functor. Then the topology introduced in \S\ref{vii.3} coincides with the smallest topology on \tannakian{\fifu} making all the elements of $\Rf'$ continuous.
\end{lemma}
\begin{proof}
The algebra of continuous functions $\Rf'$ separates points because of the existence of enough \fifu-invariant metrics. Then, for every open subset $\Omega$ with compact closure the involutive subalgebra ${\Rf'}_{\overline\Omega} \subset \C^0(\overline\Omega)$ formed by the restrictions to the closure $\overline\Omega$ of elements of $\Rf'$ is dense in the subspace $\Rf_{\overline\Omega} = \{r|_{\overline\Omega}: r \in \Rf\}$ with respect to the sup-norm, as a consequence of the \mbox{Stone}--\mbox{Weierstrass} theorem.

The subsets of the form $\tannakian{\fifu}|_{U\times U'}$, where $U$ and $U'$ are open subsets of $M$ with compact closure, are certainly open and of compact closure, as well as open relative to the topology associated with $\Rf'$. Let $\Omega$ be any one of these open subsets. We claim that both topologies agree on $\Omega$. Indeed, for each $r \in \Rf$ the restriction $r|_\Omega$ must be a uniform limit of functions which are continuous for the $\Rf'$ topology, and hence $r$ itself must be a continuous function for the $\Rf'$ topology.
\end{proof}

\begin{paragrafo}[Remark.]\label{vii.5}
Each arrow $\lambda$ of the groupoid \tannakian{\fifu} acts as a unitary transformation with respect to all \fifu-invariant metrics. More explicitly, for every $R \in \Ob(\Kt)$ and \fifu-invariant metric $\phi$ on $\fifu(R)$ one has
\begin{equazione}\label{N.v16}
\smash{\bigsca{\lambda(R) v}{\lambda(R) v'}}_\phi = \smash{\scalare{v}{v'}}_\phi\text.
\end{equazione}
\end{paragrafo}

\noindent We use this remark in the proof of the following

\begin{proposizione}\label{N.v17}
Let \fifu\ be a proper fibre functor. Then \tannakian{\fifu} is a topological groupoid.
\end{proposizione}
\begin{proof}
(a) Continuity of the inverse map $\imap$. By Lemma \ref{N.v15}, it suffices to prove that the composite ${r\circ\imap}$ is continuous for every $r = r_{R,\phi,\zeta,\zeta'}$ with $\phi$ \fifu-invariant. This is clear, for by Remark \ref{vii.5}
\begin{displaymath}
{r_{R,\phi,\zeta,\zeta'} \circ \imap} = \overline{(r_{R,\phi,\zeta',\zeta})}\text.
\end{displaymath}

(b) Continuity of the composition map \c. We start with a preliminary observation.

Let $R \in \Ob(\Kt)$. Let $\phi$ be any \fifu-invariant metric on $\fifu(R)$. For any given arrow $\lambda: x \to x'$ in \tannakian{\fifu} we can fix a local $\phi$-orthonormal frame $\zeta'_1, \ldots, \zeta'_d$ of sections defined over some neighbourhood $U'$ of $x'$. [See \S\ref{N.iiiA4}.] Choose an open neighbourhood $\Omega$ of $\lambda$ such that $\t(\Omega) \subset U'$. Let $\zeta$ be a global section of $\fifu(R)$ and let $\Phi_i$ ($i = 1, \ldots, d$) be arbitrary continuous functions on $\Omega$. The function
\begin{equazione}\label{N.v19}
\mu\: \mapsto\: \modulo{{\mu(R) \cdot \zeta(\s[\mu])} - \txtsum{i=1}d{\Phi_i(\mu) \zeta'_i(\t[\mu])}}
\end{equazione}
is certainly continuous; indeed, by (\ref{N.v16}), its square is
\begin{displaymath}
\bigmod{\zeta(\s[\mu])}^2 - 2\displaysum{i}{}{\parteRe{\left[\overline{\Phi_i(\mu)} \bigsca{\mu(R) \zeta(\s[\mu])}{\zeta'_i(\t[\mu])}\right]}} + \txtsum{i=1}d{\bigmod{\Phi_i(\mu)}^2}\text.
\end{displaymath}
Upon making the substitution $\Phi_i(\mu) = \bigsca{\mu(R) \zeta(\s[\mu])}{\zeta'_i(\t[\mu])}$ in (\ref{N.v19}), we get a function vanishing at $\lambda$ since by construction the vectors $\zeta'_i(x')$ constitute an orthonormal basis.

Now, we have to check the continuity of all functions of the form
\begin{equazione}\label{N.v21}
(\mu',\mu)\: \mapsto\: ({r_{R,\phi,\zeta,\eta} \circ \c})(\mu',\mu) = \smash{\bigsca{\mu'(R) \cdot \mu(R) \cdot \zeta(\s[\mu])}{\eta(\t[\mu'])}}_\phi
\end{equazione}
with $\phi$ \fifu-invariant. Let $x \xto\lambda x' \xto{\lambda'} x''$ be any pair of composable arrows. By the foregoing observation and \refequ{N.v16}, we see that for each $\epsilon>0$ there is a neighbourhood $\Omega_\epsilon$ of $\lambda$ such that for all composable $(\mu',\mu)$ with $\mu \in \Omega_\epsilon$ the value of the function (\ref{N.v21}) at $(\mu',\mu)$ differs from
$$%
\txtsum{i=1}d{r_{R,\phi,\zeta,\zeta'_i}(\mu) \smash{\bigsca{\mu'(R) \cdot \zeta'_i(\s[\mu'])}{\eta(\t[\mu'])}}_\phi} = \txtsum{i=1}d{r_{R,\phi,\zeta,\zeta'_i}(\mu) r_{R,\phi,\zeta'_i,\eta}(\mu')}
$$%
by $C\epsilon$ at most, where $C$ is a positive bound for the $\phi$-norm of the section $\eta$ in a given neighbourhood of $x''$.
\end{proof}

\sezione{Proof of the Reconstruction Theorem}\label{viii}
We start with some results which hold for an arbitrary Euclidean stack \stack F. We introduce the shorthand \tannakian{\G} for the Tannakian groupoid associated with the standard forgetful functor (of type \stack F) of a Lie groupoid \G.

\begin{paragrafo}[The enveloping homomorphism.]\label{viii.1}
The canonical homomorphism
\begin{equazione}\label{O.esp86}
\envelopemorph{\G}: \G \longto \tannakian{\G}
\end{equazione}
is defined by means of the identity ${\envelopemorph{\G}(g)}(E,\varrho) = \varrho(g)$. [Recall Example \ref{N.v13}.] We shall refer to \envelopemorph{\G} as the \textit{enveloping homomorphism} (\textit{of type \stack F}) \textit{of \G}.
\end{paragrafo}

\begin{theorem}\label{O.thm3}
Let \G\ be a proper Lie groupoid. Then the enveloping homomorphism of \G\ is a surjection.
\end{theorem}
\begin{proof}
To begin with, we prove that whenever $\G(x,x')$ is empty, so must be $\tannakian{\G}(x,x')$. Let $\varphi: {{\G x} \cup {\G x'}} \to \nC$ be the function which takes the value one on the orbit ${\G x}$ and the value zero on the orbit ${\G x'}$. This function is well-defined, because $\G(x,x')$ is empty. By Corollary \ref{iv8}, there is a global invariant smooth function $\Phi$ extending $\varphi$. Being invariant, $\Phi$ determines an endomorphism $a$ of the trivial representation $\TU \in {\Ob\, \R{\G}}$ such that $a_z = {\Phi(z) \id}$ for all $z$ (thus, in particular, $a_x = \id$ and $a_{x'} = 0$). Now, suppose $\lambda \in \tannakian{\G}(x,x')$. Because of the naturality of $\lambda$, the existence of the morphism $a$ contradicts the invertibility of the linear map $\lambda(\TU)$.

We are therefore reduced to proving that the induced isotropy homomorphisms $\envelopemorph{\G}|_x: \G|_x \to \tannakian{\G}|_x$ are surjective. This is now a direct consequence of Propositions \ref{O.prp5} and \ref{N.iv5}.
\end{proof}

\begin{definizione}\label{O.dfn3}
We say that a Lie groupoid \G\ is \textit{reflexive} or \textit{self-dual} (\textit{relative to \stack F}) when its enveloping homomorphism is an isomorphism of topological groupoids.
\end{definizione}

\begin{theorem}\label{O.thm1}
Let \G\ be a proper Lie groupoid. Then in order that \G\ may be reflexive it is enough that its enveloping homomorphism be injective.
\end{theorem}
\begin{proof}
The continuity of \envelopemorph{\G} is obvious, hence what we really have to show is that for each open subset $\Gamma$ of \G\ and for each point $g_0 \in \Gamma$ the image $\envelopemorph{\G}(\Gamma)$ is a neighbourhood of $\envelopemorph{\G}(g_0)$ in \tannakian{\G}.

Let $g_0 \in \G(x_0,{x_0}')$. We start by observing that it is possible to find a representation $R = (E,\varrho)$ whose associated $x_0$-th isotropy homomorphism $\varrho_0: \G|_0 \to \GL(E_0)$ is injective. [Compare the proof of Proposition \ref{O.prp5} and also \S\ref{iv.1}.] Fix an arbitrary metric $\phi$ on $E$ and local $\phi$-orthonormal frames
\begin{displaymath}
\zeta_1, \ldots, \zeta_d \quad \text{about $x_0$} \quad \text{and} \quad \zeta'_1, \ldots, \zeta'_d \quad \text{about ${x_0}'$.}
\end{displaymath}
Choose any compactly supported smooth function $0\leqq a\leqq 1$, resp.\ $0\leqq a'\leqq 1$ with support lying close enough to $x_0$, resp.\ ${x_0}'$ and such that $a(z) = 1 \aeq z = x_0$, resp.\ $a'(z) = 1 \aeq z = {x_0}'$. Then put
\begin{displaymath}
\varrho_{i,i'} \bydef {r_{i,i'} \circ \envelopemorph{\G}} \bydef {r_{R,\phi,\zeta_i,\zeta'_{\smash{i'}}} \circ \envelopemorph{\G}}\text, \quad \text{and} \quad \varrho_{\iota,\iota'} \bydef {r_{\iota,\iota'} \circ \envelopemorph{\G}} \bydef {({a \circ \s_{\G}})({a' \circ \t_{\G}})}
\end{displaymath}
with $\iota=0$ or $\iota'=0$. Finally, let $\omega_{\iota,\iota'} = \varrho_{\iota,\iota'}(g_0)$ for $0\leqq \iota,\iota'\leqq d$.

We claim that there exist open disks $D_{\iota,\iota'}$, with $D_{\iota,\iota'}$ encircling the complex number $\omega_{\iota,\iota'}$, which satisfy
\begin{equazione}\label{xx6}
\displaycap{0\leqq \iota,\iota'\leqq d}{}{{\varrho_{\iota,\iota'}}^{-1}(D_{\iota,\iota'})} \subset \Gamma\text.
\end{equazione}
Once this claim is proven, the statement that $\envelopemorph{\G}(\Gamma)$ is a neighbourhood of $\envelopemorph{\G}(g_0)$ will be proven as well. Indeed, by Theorem \refcnt{O.thm3} we have
\begin{displaymath}
{\bigcap {r_{\iota,\iota'}}^{-1}(D_{\iota,\iota'})}\: =\: {\envelopemorph\G\, {\envelopemorph\G}^{-1} \left(\smash[b]{\displaycap{}{}{{r_{\iota,\iota'}}^{-1} (D_{\iota,\iota'})}}\right)}\: =\: {\envelopemorph\G \left(\smash[b]{\displaycap{}{}{{\varrho_{\iota,\iota'}}^{-1} (D_{\iota,\iota'})}}\right)}
\end{displaymath}
where each ${r_{\iota,\iota'}}^{-1}(D_{\iota,\iota'})$ is an open neighbourhood of $\envelopemorph{\G}(g_0)$ in \tannakian{\G}.

In order to establish (\ref{xx6}), we fix for each $0\leqq \iota,\iota'\leqq d$ a decreasing sequence of open disks centred at $\omega_{\iota,\iota'}$
\begin{equazione}\label{xx7}
\cdots \subset {D_{\iota,\iota'}}^{p+1} \subset {D_{\iota,\iota'}}^p \subset \cdots \subset {D_{\iota,\iota'}}^1 \subset \nC
\end{equazione}
with radius converging to zero. If we agree that ${D_{\iota,\iota'}}^1$ has radius $\tfrac12$ then
\begin{equazione+}\label{xx8}
\Sigma^p\: \bydef\: {\bigcap {{r_{\iota,\iota'}}^{-1} \bigl(\overline{{D_{\iota,\iota'}}^p}\bigr)} - \Gamma} & (p = 1, 2, \ldots)
\end{equazione+}
is a closed subset of the compact space $\G(K,K')$ where $K = \support a$ and $K' = \support{a'}$. The intersection \txtcap{p=1}{\infty}{\Sigma^p} is empty because of the injectivity of the map $\G(x_0,{x_0}') \to \Iso(E_{x_0},E_{{x_0}'})\text,\: g \mapsto \varrho(g)$ and the choice of $a,a'$. Thus, there must be some $p$ such that $\Sigma^p = \varnothing$. This proves the claim.
\end{proof}

\noindent In \S\ref{N.v13}, we remarked on passing that \envelopemorph{\G} is a morphism of $\C^\infty$-spaces. This is in fact true for an arbitrary, not necessarily proper Lie groupoid \G. One may wonder whether more can be said when \G\ is reflexive.

Hereafter we shall freely make use of some notation introduced in the context of the preceding proof. We define the smooth mappings
\begin{multiriga}\label{xx9}
\varrho^{\zeta_1,\ldots,\zeta_d}_{\zeta'_1,\ldots,\zeta'_d}: \G \longto {M \times M \times \End(\nC^d)}\text,\\
g \mapsto \bigl(\s(g);\t(g);\varrho_{1,1}(g),\ldots,\varrho_{i,i'}(g),\ldots,\varrho_{d,d}(g)\bigr)\text,
\end{multiriga}
where $M$ is the base of \G, and introduce the abbreviations $\zeta \equiv \zeta_1,\ldots,\zeta_d$, $\zeta \equiv \zeta'_1,\ldots,\zeta'_d$. If the homomorphism \envelopemorph{\G} is faithful, Lemma \ref{O.lem8} implies that for each arrow $g_0$ there exists a representation $R = (E,\varrho)$ such that the map $\G(x_0,{x_0}') \longto \Iso(E_{x_0},E_{{x_0}'})\text,\: g \mapsto \varrho(g)$ becomes injective when restricted to a sufficiently small open neighbourhood of $g_0$.

\begin{lemma}\label{O.lem1}
Suppose the map $\G(x_0,{x_0}') \to \Iso(E_{x_0},E_{{x_0}'})\text,\: g \mapsto \varrho(g)$ is injective near $g_0$. Then \refequ{xx9} is an immersion at $g_0$.
\end{lemma}
\begin{proof}
Fix open balls $U$ and $U'$ centred at $x_0$ and ${x_0}'$ respectively, so small that the sections $\zeta_1, \ldots, \zeta_d$ (resp.\ $\zeta_1', \ldots, \zeta_d'$) form a local orthonormal frame for $E$ over $U$ (resp.\ $U'$). Up to a local diffeomorphism, the map \refequ{xx9} has the following form near $g_0$, provided $U$ is chosen small enough:
\begin{equazione}\label{xx.10}
{U \times \nR^k} \to {U \times U' \times \End(\nC^d)}\text, \quad (u,v) \mapsto \bigl(u;u'(u,v);\boldsymbol\varrho(u,v)\bigr)\text,
\end{equazione}
where $\boldsymbol\varrho(g)$ denotes the matrix $\{\varrho_{i,i'}(g)\}_{1\leqq i,i'\leqq d}$. Evidently, \refequ{xx.10} is immersive at $g_0 = (x_0,0)$ if and only if the partial map $v \mapsto \bigl(u'(x_0,v);\boldsymbol\varrho(x_0,v)\bigr)$ is immersive at zero. We are therefore reduced to showing that the restriction of \refequ{xx9} to $\G(x_0,\text-)$ is immersive at $g_0$.

Let $G$ be the isotropy group of \G\ at $x_0$. By choosing a local equivariant trivialization $\G(x_0,S) \iso {S \times G}$ where $S$ is a submanifold of $U'$ passing through ${x_0}'$, the restriction of \refequ{xx9} to $\G(x_0,\text-)$ takes the form
\begin{equazione}\label{xx.12}
{S \times G} \to {U' \times \End(\nC^d)}\text, \quad (s,g) \mapsto \bigl(s;\boldsymbol\varrho(s,g)\bigr)\text.
\end{equazione}
This map is immersive at $g_0 = ({x_0}',e)$ if and only if so is at $e$ the partial map $g \mapsto \boldsymbol\varrho({x_0}',g)$, where $e$ is the unit of the group $G$. Thus, it suffices to show that the isotropy representation $G \to \GL(E_{x_0})$ induced by $\varrho$ is immersive at $e$. By hypothesis, this representation is injective in an open neighbourhood of $e$ and hence our claim follows at once.
\end{proof}

\noindent Let an arrow $\lambda_0 \in \tannakian{\G}$ be given. We contend that there exists some open neighbourhood $\Omega$ of $\lambda_0$ such that $(\Omega,\Rf^\infty_\Omega)$ is isomorphic, as a $\C^\infty$-space, to a smooth manifold $(X,\smooth[X])$.

Since \G\ is reflexive, there is a unique $g_0 \in \G$ such that $\lambda_0 = \envelopemorph{\G}(g_0)$. By Lemma \refcnt{O.lem1} and the remarks preceding it, we can find some $R$ for which there exists an open neighbourhood $\Gamma$ of $g_0$ in \G\ such that $\varrho^\zeta_{\zeta'}$ induces a diffeomorphism of $\Gamma$ onto a submanifold $X$ of ${M \times M \times \End(\nC^d)}$. Define
\begin{multiriga}\label{xx.13}
r^{\zeta_1\ldots,\zeta_d}_{\zeta_1',\ldots,\zeta_d'}: \tannakian{\G} \longto {M \times M \times \End(\nC^d)}\text,\\
\lambda \mapsto \bigl(\s(\lambda);\t(\lambda);r_{1,1}(\lambda),\ldots,r_{i,i'}(\lambda),\ldots,r_{d,d}(\lambda)\bigr)\text.
\end{multiriga}
This map is evidently a morphism of $\C^\infty$-spaces. By the reflexivity of \G, \envelopemorph{\G} induces a homeomorphism between $\Gamma$ and the open subset $\Omega \equiv \envelopemorph{\G}(\Gamma)$ of \tannakian{\G}. Clearly, ${\varrho^\zeta_{\zeta'}|_\Gamma} = {{r^\zeta_{\zeta'}|_\Omega} \circ {\envelopemorph{\G}|_\Gamma}}$ and so ${r^\zeta_{\zeta'}|_\Omega}$ yields a homeomorphism between $\Omega$ and $X$.

We claim that the map $r^\zeta_{\zeta'}|_\Omega$ is the desired isomorphism of $\C^\infty$-spaces. \textsc{(a)} In one direction, suppose $f \in \C^\infty(X)$. Because of the local character of the claim, it is no loss of generality to assume that $f$ admits a smooth extension
$$%
\tilde f \in \C^\infty\bigl({M \times M \times \End(\nC^d)}\bigr)\text,
$$%
thus ${f \circ {r^\zeta_{\zeta'}|_\Omega}} = {\tilde f \circ {r^\zeta_{\zeta'}|_\Omega}}$ is evidently an element of $\Rf^\infty(\Omega)$. \textsc{(b)} Conversely, let $f: X \to \nC$ be a function such that ${f \circ {r^\zeta_{\zeta'}|_\Omega}}$ belongs to $\Rf^\infty(\Omega)$. Since \envelopemorph{\G} is a morphism of $\C^\infty$-spaces, the composite ${f \circ {r^\zeta_{\zeta'}|_\Omega} \circ {\envelopemorph{\G}|_\Gamma}} = {f \circ {\varrho^\zeta_{\zeta'}|_\Gamma}}$ will belong to $\C^\infty(\Gamma)$. As ${\varrho^\zeta_{\zeta'}|_\Gamma}$ is a diffeomorphism, it follows that $f \in \C^\infty(X)$. The claim is proven.

Summarizing our conclusions:

\begin{proposizione}\label{O.prp10}
Let \G\ be a reflexive groupoid (Definition \ref{O.dfn3}). Then the enveloping homomorphism \envelopemorph{\G} is an isomorphism of $\C^\infty$-spaces; it follows that the Tannakian groupoid \tannakian{\G} is a Hausdorff Lie groupoid, isomorphic to \G.
\end{proposizione}

\noindent We shall now turn our attention to a very delicate issue, namely the injectivity of the enveloping homomorphism. Clearly, \envelopemorph{\G} is injective if and only if \G\ admits enough representations; this means that for each $x \in M$ and $g \neq x$ in the $x$-th isotropy group of \G\ there is a representation $(E,\varrho)$ such that $\varrho(g) \neq \id \in \Aut(E_x)$. For a generic Lie groupoid \G, this property dramatically depends on the type of representations one is considering.

We claim that each proper Lie groupoid admits enough representations on smooth Euclidean fields (cfr \S\ref{xvi.2}). For the rest of the section, we shall exclusively deal with such representations.

\begin{paragrafo}[Cut-off functions.]\label{viii.7}
We begin with some preliminary remarks of a purely topological nature. Let \G\ be a proper Lie groupoid over a manifold $M$. Recall that a subset $S \subset M$ is said to be invariant when ${s \in S} \seq {{g \cdot s} \in S}$ for all arrows $g$. If $S$ is any subset of $M$, we let ${\G \cdot S}$ denote the saturation of $S$, that is to say the smallest invariant subset of $M$ containing $S$. The saturation of an open subset is also open. It is an easy exercise to show that ${\G \cdot \overline V} = \overline{\G \cdot V}$ for all open subsets $V$ with compact closure. It follows that if $U$ is an invariant open subset of $M$ then $U$ coincides with the union over all invariant open subsets $V$ whose closure is compact and contained in $U$. The last remark applies to the construction of \G-invariant partitions of unity over $M$; for our purposes, it will be enough to illustrate a special case of this construction. Consider an arbitrary point $x_0 \in M$ and let $U$ be an open invariant neighbourhood of $x_0$. Choose another open neighbourhood $V$ of $x_0$, invariant and with closure contained in $U$. The orbit ${\G \cdot x_0}$ and the set-theoretic complement ${\complement V}$ are invariant disjoint closed subsets of $M$, so by Corollary \ref{iv8} there exists an invariant smooth function on $M$ which takes the value one at $x_0$ and vanishes outside $V$.
\end{paragrafo}

\begin{paragrafo}[Extendability of proper Lie groupoid actions on smooth Euclidean fields.]\label{viii.8}
Let \G\ be a proper Lie groupoid, with base $M$. Suppose we are given a ``partial'' representation $(\E_U,\varrho_U)$ of $\G|_U$ on a smooth Euclidean field $\E_U$ over $U$, where $U$ is an invariant open neighbourhood of a point $x_0$ in $M$. We want to show that there exists a ``global'' representation $(\E,\varrho)$ of \G\ on a smooth Euclidean field \E\ such that $(\E_U)_0 \equiv (\E_U)_{x_0}$ and $\E_0 \equiv \E_{x_0}$ are isomorphic $G$-modules, where $G$ is the isotropy group of \G\ at $x_0$.

To begin with, we fix any invariant smooth function $a \in \C^\infty(M)$ with $a(x_0) = 1$ and $\support a \subset U$ (cut-off function). Let $V$ denote the set of all $x$ such that $a(x)\neq 0$. Define $\E_x$ to be the fibre $(\E_U)_x$ if $x \in V$ and $\{0\}$ otherwise. Let \f{E} be the following sheaf of sections of the bundle $\{\E_x\}$:
\begin{equazione}\label{N.vi4}
W\: \mapsto\: \bigl\{\text{``prolongation of ${a \zeta}$ by zero''}: \zeta \in \f{(E_{\mathnormal U})}({U \cap W})\bigr\}\text.
\end{equazione}
These data define a smooth Euclidean field \E\ over $M$. Define $\varrho(g)$ to be $\varrho_U(g)$ if $g \in \G|_V$ and the zero map otherwise. The bundle of linear maps
$$%
\bigl\{\varrho(g): ({\s^*\E})_g \isoto ({\t^*\E})_g\bigr\}
$$%
will provide an action of \G\ on \E\ as long as it is a morphism of smooth Euclidean fields over \G\ of ${\s^*\E}$ into ${\t^*\E}$. Now, by the invariance of $a$ and the local expression (\ref{N.iv1c}) for $\varrho_U$, one has
$$%
{\varrho(g) [{a \zeta}(\s[g])]} = {a(\s[g]) \varrho(g) \zeta(\s[g])} = {a(\t[g]) \txtsum{i=1}{d}{r_i(g) \zeta'_i(\t[g])}} = \txtsum{i=1}{d}{r_i(g) [{a \zeta'_i}(\t[g])]}\text,
$$%
as desired. Finally, the identity $\E_0 = (\E_U)_{x_0}$ (by construction) is a $G$-equivariant isomorphism.
\end{paragrafo}

\noindent Putting Theorem \ref{O.thm1}, Proposition \ref{O.prp10}, the considerations of \S\ref{ZungThm} and those of the last subsection together, we conclude

\begin{theorem}[Reconstruction Theorem]\label{xx.18}
Within the type $\Euc^\infty$ of smooth Euclidean fields, every proper Lie groupoid is reflexive, that is to say $\C^\infty$-isomorphic to its Tannakian groupoid via the corresponding enveloping homomorphism.
\end{theorem}
\begin{proof}
A faithful representation $G \into \GL(\boldsymbol E)$ of a compact Lie group $G$ on a finite dimensional vector space $\boldsymbol E$ induces, for any smooth action of $G$ on a smooth manifold $V$, a faithful representation of the action groupoid ${G \ltimes V}$ on the trivial vector bundle ${V \times \boldsymbol E}$.
\end{proof}

\bibliography{../../SharedTeX/Biblio/groupoid,../../SharedTeX/Biblio/tannaka}

\end{document}